\documentclass[twoside]{book}
\usepackage[top=0.9in, bottom=0.9in,left=0.9in, right=0.9in, paperwidth=6in, paperheight=9in]{geometry}

\usepackage[russian,english]{babel}
\usepackage[utf8]{inputenc}

\usepackage{euclid2alexandrov-pdp}
\hypersetup{pdftitle={Piecewise distance preserving maps.
},%
pdfauthor={Anton Petrunin and Allan Yashinski}}

\usepackage{fancyhdr}
\addto\captionsenglish{%
 }
 
\usepackage[nottoc]{tocbibind}

\makeindex
\begin{document}
\title{Lectures on\\piecewise distance-preserving maps}
\author{Anton Petrunin and Allan Yashinski}
\date{}
\maketitle
\chapter*{Preface}
These lectures were part of a geometry course held during the Fall 2011 Mathematics Advanced Study Semesters (MASS) Program at Penn State.

The lectures are meant to be accessible to advanced undergraduate and beginning graduate students in mathematics.
We have placed a great emphasis on clarity and exposition, and we have included many exercises.
Hints and solutions for most of the exercises are provided at the end.

The lectures discuss piecewise distance-preserving maps from a 2-dimensional polyhedral space into the plane.
Informally speaking, a polyhedral space is a space that is glued together out of triangles, for example the surface of a polyhedron.
If one imagines such a polyhedral space as a paper model, then a piecewise distance-preserving map into the plane is essentially a way to fold the model so that it lies flat on a table.
We have five lectures on the following topics:

\begin{itemize}
\item Zalgaller's folding theorem, which guarantees the existence of a piecewise distance-preserving map from a 2-dimensional polyhedral space into the plane.
In other words, it is always possible to fold the paper model onto the table.
\item Brehm's extension theorem, which allows one to build piecewise distance-preserving maps from a convex polygon into the plane with prescribed images on a finite subset of the polygon.
\item Akopyan's approximation theorem, which allows one to approximate maps from a 2-dimensional polyhedral space into the plane by piecewise distance-preserving maps.
\item Gromov's rumpling theorem, which shows the existence of a length-preserving map from the sphere into the plane,
i.e., a map that preserves the lengths of all curves.

{\sloppy

\item An entertaining problem posed by Arnold on paper folding, which asks whether it is possible to fold a square in the plane so that its perimeter increases.

}
\end{itemize}

We consider only the 2-dimensional case to keep things easy to visualize.
However, 
most of the results admit generalizations to higher dimensions.
These results are discussed in the final remarks, where proper credit and references are given.

{\sloppy

\parbf{Acknowledgments.} 
We would like to thank
Arseniy Akopyan, 
Robert Lang,
Vansh Sehrawat,
and Alexei Tarasov
for their help.
We would also like to thank all the students in our class
for their participation and true interest.

}

\tableofcontents

\chapter*{Preliminaries}
\addcontentsline{toc}{chapter}{Preliminaries}
\addtocontents{toc}{\protect\begin{quote}}

This chapter serves as a quick review of the necessary background material.
Its primary purpose is to refresh the reader's memory and to familiarize the reader with our notation.
For a more in-depth account of this material, we recommend the first three chapters of the book
``Metric Geometry'' by Burago--Burago--Ivanov \cite{BBI}.

\section*{Metric spaces}
\addtocontents{toc}{Metric spaces.}

\begin{thm}{Definition}
\label{def:metric-space}
A metric space is a pair $(X,\Dist)$ where $X$ is a set and 
$$\Dist \: X \times X \rightarrow [0,\infty)$$
is a function satisfying

\begin{subthm}{def:metric-space:zero}
$\Dist(x,y) = 0$ if and only if $x=y$,
\end{subthm}

\begin{subthm}{}
$\Dist(x,y) = \Dist(y,x)$ for any $x, y \in X$,
\end{subthm}

\begin{subthm}{triangle-inq}
$\Dist(x,z)\le \Dist(x,y)+\Dist(y,z)$ for any $x, y, z \in X$.
\end{subthm}

\end{thm}

Condition (\ref{SHORT.triangle-inq}) is called the \emph{triangle inequality}\index{triangle inequality}.
An element $x \in X$ is called a \emph{point}\index{point} of the metric space $X$.
The function $\Dist\: X \times X \z\rightarrow [0,\infty)$ is called a \emph{metric}\index{metric},
and the non-negative value $\Dist(x,y)$ is called the \emph{distance}\index{distance} from $x$ to $y$.

Most of the time, it will be understood which metric we are using on a given set.
In such cases, we may refer to the \emph{metric space $X$}, instead of the \emph{metric space $(X,\Dist)$}.
We will often use the notation $|x-y|$ to denote%
\footnote{Be aware that this is just alternate notation for $\Dist(x,y)$, and $x - y$ itself has no meaning in a general metric space.}
the distance $\Dist(x,y)$, and we shall occasionally write $|x-y|_X$ to emphasize that we are using the metric on the space~$X$.

An important example for us is $n$-dimensional Euclidean space $\RR^n$ with the standard metric
\[|x-y|_{\RR^n} = \sqrt{(x_1 - y_1)^2+\dots+(x_n - y_n)^2},\]
defined for the points $x = (x_1, \ldots , x_n)$ and $y = (y_1, \ldots , y_n)$ in $\RR^n$.

Any subset of a metric space is also a metric space, by restricting the original metric to the subset.
In this way, all subsets of Euclidean space, in particular convex polyhedra,
are metric spaces.
The unit sphere $\SS^n$, which is the set of all unit vectors in $\RR^{n+1}$, is also a metric space in this way.
However, we shall also be interested in a different metric on $\SS^n$, the \emph{angle metric}\index{angle metric}.\label{angle-metric}
This metric can be formally defined as
\[|x - y|_{\SS^n} = \arccos\langle x, y\rangle,\]
where $\langle x, y \rangle$ is the standard inner product of the two unit vectors $x, y \in \RR^{n+1}$.
The relationship between these two metrics on $\SS^n$ will be described below, in the discussion of induced length metrics.

\begin{thm}{Definition}
Let $X$ and $Y$ be metric spaces.
\begin{subthm}{}
A map $f\:X\to Y$ is \emph{distance non-expanding}\index{distance non-expanding map} if
$$|f(x)-f(x')|_Y\le |x-x'|_X$$
for any $x,x'\in X$.
\end{subthm}

\begin{subthm}{}
A map $f\:X\to Y$ is \emph{distance-preserving}\index{distance-preserving map} if
$$|f(x)-f(x')|_Y= |x-x'|_X$$
for any $x,x'\in X$.
\end{subthm}

\begin{subthm}{}
A distance-preserving bijection $f\:X\to Y$ is called an \emph{isometry}\index{isometry}.
\end{subthm}

\begin{subthm}{}
The spaces $X$ and $Y$ are called \emph{isometric}\index{isometric spaces} (briefly $X\iso Y$)
 if there is an isometry $f\:X\to Y$.
\end{subthm}

\end{thm}

\begin{thm}{Exercise}\label{ex:IsometriesOfR2}

\begin{subthm}{} \label{ex:IsometriesOfR2Uniqueness}
Prove that an isometry $\iota: \RR^2 \to \RR^2$ is uniquely determined by its effect on any three non-collinear points.
\end{subthm}

\begin{subthm}{}\label{ex:IsometriesOfR2Existence}
Show that if $x, y, z$ and $x', y', z'$ are two triples of non-collinear points in $\RR^2$ satisfying
$$|x - y| = |x' - y'|, \quad |x - z| = |x' - z'|, \quad |y - z| = |y' - z'|,$$
then there is an isometry $\iota$ of $\RR^2$ such that
\[ \iota(x) = x', \quad \iota(y) = y', \quad \iota(z) = z'.\]
By part a), $\iota$ is unique.
\end{subthm}
\end{thm}

\section*{Calculus}
\addtocontents{toc}{Calculus.}

\begin{thm}{Definition}
 Let $X$ be a metric space.
A sequence of points $x_1, x_2, \ldots$ in $X$ is called \emph{convergent}\index{convergent}
if there is 
$x_\infty\in X$ such that $|x_\infty -x_n|\to 0$ as $n\to\infty$.
That is, for every $\eps > 0$, there is a natural number $N$ such that for all $n \ge N$, we have $$|x_\infty-x_n| < \eps.$$

In this case, we say that the sequence $(x_n)$ \emph{converges} to $x_\infty$,
or $x_\infty$ is the \emph{limit} of the sequence $(x_n)$.
Notationally, we write $x_n\to x_\infty$ as $n\to\infty$
or $x_\infty=\lim_{n\to\infty} x_n$.
\end{thm}

\begin{thm}{Definition}
Let $X$ and $Y$ be metric spaces.
A map $f\:X\to Y$ is called continuous if for any convergent sequence $x_n\to x_\infty$ in $X$,
%the sequence $y_n\z=f(x_n)$ converges to $y_\infty=f(x_\infty)$ in $Y$.
we have $f(x_n) \to f(x_\infty)$ in $Y$.

Equivalently, $f\:X\to Y$ is continuous if for any $x\in X$ and any $\eps>0$,
there is $\delta>0$ such that 
$$|x-x'|_X<\delta\ \text{ implies }\ |f(x)-f(x')|_Y<\eps.$$

\end{thm}

\begin{thm}{Definition}
Let $X$ and $Y$ be metric spaces.
A continuous bijection $f\:X\to Y$ 
is called a \emph{homeomorphism}\index{homeomorphism} 
if its inverse $f^{-1}\:Y\z\to X$ is also continuous.

If there exists a homeomorphism $f\:X\to Y$,
we say that $X$ is \emph{homeomorphic}\index{homeomorphic} to $Y$,
or $X$ and $Y$ are \emph{homeomorphic}.
\end{thm}
Notice that a distance non-expanding map is always continuous, and an isometry is an example of a homeomorphism.

\begin{thm}{Definition}
A subset $A$ of a metric space $X$ is called \emph{closed}\index{closed set} if whenever a sequence $(x_n)$ of points from $A$ converges in $X$, we have that $\lim_{n\to\infty} x_n \in A$.

A set $\Omega \subset X$ is called \emph{open}\index{open set} if the complement $X \setminus \Omega$ is a closed set.
Equivalently, $\Omega \subset X$ is open if for any $z\in \Omega$, 
there is $\eps>0$ such that if $|x-z|<\eps$, then $x \in \Omega$.
\end{thm}

Note that the intersection of an arbitrary family of closed sets is closed.
It follows that for any set $Q$ in a metric space $X$, there is a minimal closed set which contains $Q$.
This set is called the \emph{closure of $Q$}\index{closure} and is denoted $\Closure Q$.
The closure of $Q$ can be obtained as the intersection of all closed sets $A$ containing $Q$.%
\footnote{Notice that the whole space $X$ is a closed subset, and hence there is at least one closed set which contains $Q$.}
The closure of $Q$ is also equal to the set of all limits of all sequences in $Q$.

Similarly, the union of an arbitrary family of open sets is open.
It follows that for any set $Q$ in a metric space, there is a maximal open set which is contained in $Q$.
This set is called the \emph{interior of $Q$}\index{interior} and is denoted by $\Interior Q$.
The interior of $Q$ can be obtained as the union of all open sets contained in $Q$.

The set-theoretic difference
$$\partial Q=\Closure Q\backslash \Interior Q$$
is called the \emph{boundary}\index{boundary} of $Q$.
Since the boundary of $Q$ depends on the space in which it is embedded, we may use the notation $\partial_X Q$ if we want to emphasize that $Q$ is a subset of the metric space $X$.
Notice that $\partial Q$ is a closed set.
One can show that a point $p$ is in $\partial Q$ if and only if for any $\eps>0$, there are points $q\in Q$ and $q'\notin Q$ such that $|p-q| < \eps$ and $|p-q'|<\eps$.

\section*{Curves}
\addtocontents{toc}{Curves.}

A \emph{real interval}\index{real interval} is a convex subset of $\RR$ which contains more than one point.
Examples include $(0,1)$, $(-\infty, 0]$, and $\RR$.

\begin{thm}{Definition}\label{def:curve}
A \emph{curve}\index{curve} is a continuous mapping $\alpha\:\II\to X$,
where $\II$ is a real interval and $X$ is a metric space.

If $\II=[a,b]$ and $$\alpha(a)=p,\ \ \alpha(b)=q,$$
we say that $\alpha$ is a \emph{curve from $p$ to $q$}\index{curve from $p$ to $q$}.
\end{thm}

\begin{thm}{Definition}\label{def:length}
Let $\alpha\:\II\to X$ be a curve.
Define the \emph{length}\index{length of curve} of $\alpha$ to be
\begin{align*}
\length \alpha
&= 
\sup \{|\alpha(t_0)-\alpha(t_1)|+|\alpha(t_1)-\alpha(t_2)|+\dots
\\
&\ \ \ \ \ \ \ \ \ \ \ \ \ \ \ \ \ \ \ \ \ \ \ \ \ \ \ \ \ \ \ \ \dots+|\alpha(t_{k-1})-\alpha(t_k)|\}, 
\end{align*}
where the supremum is taken over all positive integers $k$ and all sequences $t_0 < t_1 < \cdots < t_k$ in $\II$.

A curve is called \emph{rectifiable}\index{rectifiable curve} if its length is finite.
\end{thm}

\begin{thm}{Semicontinuity of length}\label{thm:length-semicont}
Length is a lower semicontinuous functional on the space of curves
$\alpha\:\II\to X$ with respect to pointwise convergence.

In other words: assume that a sequence
of curves $\alpha_n\:\II\to X$ converges pointwise
to a curve $\alpha_\infty\:\II\to X$;
i.e., for any fixed $t\in\II$, we have $\alpha_n(t)\to\alpha_\infty(t)$ as $n\to\infty$.
Then 
$$\liminf_{n\to\infty} \length\alpha_n \ge \length\alpha_\infty.\eqlbl{eq:semicont-length}$$

\end{thm}

See \cite[Proposition 2.3.4]{BBI} for a proof.

\begin{wrapfigure}{r}{32mm}
\begin{lpic}[t(-8mm),b(-3mm),r(0mm),l(0mm)]{pics/stairs(0.5)}
%\lbl[lb]{162,180;$\alpha_0$}
%\lbl[lb]{82,150;{\small $\alpha_1$}}
%\lbl[lb]{42,130;{\tiny $\alpha_2$}}
\end{lpic}
\end{wrapfigure}

Note that the inequality \ref{eq:semicont-length} might be strict.
For example,
the diagonal $\alpha_\infty$ of the unit square (gray in the picture)
can be approximated by a sequence of staircase-like polygonal curves $\alpha_n$
with sides parallel to the sides of the square (for example, $\alpha_6$ is black in the picture).
In this case
$\length\alpha_\infty=\sqrt{2}$
and $\length\alpha_n=2$ for all $n$.

%Note that applying triangle inequality to the definition of length, 
%we get that if $\gamma$ is a curve from $x$ to $y$ in a metric space $X$ 
%then
%$$\length \gamma\ge |x-y|_X.$$

By taking $t_0$ and $t_1$ to be the endpoints of $\II$ in the definition of length, it follows that
$$\length \alpha\ge |x-y|_X$$ whenever $\alpha$ is a curve in $X$ from $x$ to $y$.

\begin{thm}{Definition}\label{def:length-space}
A metric space $X$ is called a \emph{length space}\index{length space} if for any two points $x,y\in X$ and any $\eps>0$, there is a curve $\alpha$ from $x$ to $y$ such that
$$\length \alpha<|x-y|_X+\eps.$$

\end{thm}

Let $X$ be a metric space.
Consider the function $\hat d\:X\z\times X\z\to\RR \cup \{\infty\}$ defined by
$$\hat d(x,y)\df\inf_\alpha\{\length\alpha\}$$
where the infimum is taken over all curves $\alpha$ from $x$ to $y$
(if there is no such curve, then $\hat d(x,y)=\infty$.)

It is straightforward to see that $\hat d\:X\times X\to\RR$ satisfies all conditions of a metric, provided $\hat d(x,y)<\infty$ for all $x,y\in X$.
In this case, the metric $\hat d$ will be called the \emph{induced length metric}\index{induced length metric} of the metric $|{*}-{*}|_X$.
By construction, $(X, \hat d)$ is a length space.

For example, the angle metric on $\SS^n$ discussed on page \pageref{angle-metric}
is the induced length metric of the restriction of the Euclidean metric on $\RR^{n+1}$ to $\SS^n$.

\begin{thm}{Definition}\label{def:length-preserving}
A continuous map $f\:X\z\to Y$ between two metric spaces
is called \emph{length-preserving}\index{length-preserving map} if for any curve $\alpha\:\II\to X$,
we have 
$$\length\alpha=\length(f\circ\alpha).$$

\end{thm}

Note that since $f$ is continuous, the composition $f\circ\alpha$ is continuous, hence a curve in $Y$.
Therefore the above definition makes sense.

\begin{thm}{Exercise}\label{LP=>short}
Let $X$ and $Y$ be length spaces.
Show that 

\begin{enumerate}[a)]
\item\label{LP=>short:a} Any length-preserving map $f\: X\to Y$
is also distance non-expanding.
\item\label{LP=>short:b} A distance non-expanding map $f\: X\to Y$ is length-preserving if
for any two points $p$ and $q$ in $X$
and any curve $\alpha$ from $p$ to $q$, we have 
$$\length(f\circ\alpha)\ge |p-q|.$$
\end{enumerate}
\end{thm}

\begin{thm}{Definition}\label{def:geodesic}
 A curve $\alpha\:\II\to X$ is called a \emph{geodesic}\index{geodesic}%
\footnote{Formally, our ``geodesic'' should be called ``unit-speed minimizing geodesic'',
and the term ``geodesic'' is reserved for curves which \emph{locally} satisfy the identity
$$|\alpha(t_0)-\alpha(t_1)|_X=\Const\cdot|t_0-t_1|$$
for some $\Const\ge 0$.}
 if it is a distance-preserving map;
i.e., if 
$$|\alpha(t_0)-\alpha(t_1)|_X=|t_0-t_1|$$
for any $t_0,t_1\in\II$.
\end{thm}

If $\alpha$ is a geodesic from $p$ to $q$, then the image $\alpha(\II)$
will also be denoted by $[p,q]$.
Once we write $[p,q]$, we mean that there is at least one geodesic from $p$ to $q$ and we made a choice of one of them.

\section*{Polyhedral spaces}
\addtocontents{toc}{Polyhedral spaces.}

A subset $C \subseteq \RR^n$ is \emph{convex} if for any two points $x, y \in C$, the line segment connecting $x$ and $y$ lies entirely in $C$.
In other words,
\[(1-t)\cdot x + t\cdot y\in C\]
for any $t\in[0,1]$.

Given a subset $V\subset\RR^n$, 
the intersection of all convex sets containing $V$ is called 
the \emph{convex hull}\index{convex hull} of $V$,
and will be denoted by $\Conv V$.

The convex hull of a finite subset of $\RR^n$ is called a \emph{convex polyhedron}\index{convex polyhedron}.
A convex polyhedron is a metric space under the metric it inherits as a subset of $\RR^n$.

\parbf{Simplices.}
Assume $V=\{v_0,\dots,v_m\}$ is a finite subset of $\RR^n$
such that the $m$ vectors 
$$v_1-v_0,\ v_2-v_0,\ \dots,\ v_m-v_0$$
are linearly independent.
Then the convex hull $\Delta^m=\Conv V$
is called an \emph{$m$-dimensional (Euclidean) simplex}\index{simplex}.

So,
a 0-dimensional simplex is a one-point set; 
a 1-dimensional simplex is a line segment;
a 2-dimensional simplex is a triangle;
a 3-dimensional simplex is a tetrahedron.

If $\Delta^m$ is as above,
then the convex hull of any $(k+1)$-point subset of $\{v_0,\dots,v_m\}$ is a $k$-dimensional simplex, which will be called a \emph{face}\index{face} of $\Delta^m$.

\parbf{Barycentric coordinates.}
Let $\Delta^m=\Conv\{v_0,\dots,v_m\}$ be an $m$-dimensional simplex.
Note that $x\in \Delta^m$ if and only if 
$$x=\lambda_0\cdot v_0+\dots+\lambda_m\cdot v_m$$
for some (necessarily unique) array of non-negative real numbers $\lambda_0,\z\dots,\lambda_m$ such that
\[\lambda_0+\dots+\lambda_m=1.\]
In this case, the real array $(\lambda_0,\dots,\lambda_m)$ will be called the \emph{barycentric coordinates}\index{barycentric coordinates} of the point $x$.

\parbf{Simplicial complexes.}
A \emph{simplicial complex}\index{simplicial complex} is defined as a finite collection $\mathcal{K}$
of simplices in $\RR^n$ that satisfies the following conditions:
\begin{itemize}
\item Any face of a simplex from $\mathcal{K}$ is also in $\mathcal{K}$.
\item The intersection of any two simplices $\Delta_1$ and $\Delta_2\in \mathcal{K}$ is either the empty set or a face of both $\Delta_1$ and $\Delta_2$.
\end{itemize}

The dimension of the simplicial complex $\mathcal{K}$
(briefly $\dim \mathcal{K}$)
is defined as the maximal dimension of all of its simplices.

For example, a $1$-dimensional simplicial complex,
also called a \emph{graph},
is a finite collection of points 
with a collection of non-crossing edges connecting some of these points.
An example of a $2$-dimensional simplicial complex is the surface of a tetrahedron in $\RR^3$, or, more generally, the surface of any polyhedron in $\RR^3$.
Two disjoint $2$-simplices with a single segment connecting a vertex from each simplex give another example of a $2$-dimensional simplicial complex.
Note that in the last example 
not every simplex is a face of a simplex of maximal dimension.

We say that a point $x$ belongs to the simplicial complex $\mathcal{K}$ if it belongs to one of its simplices.
The set of all points of $\mathcal{K}$ is called the \emph{underlying set} of $\mathcal{K}$,
which will be denoted by $|\mathcal{K}|$.
Since $|\mathcal{K}| \subseteq \RR^n$, $|\mathcal{K}|$ is naturally a metric space.

A metric space $X$ is called a \emph{topological polytope}
if there is a simplicial complex $\mathcal{K}$
and a homeomorphism
$f\:|\mathcal{K}|\to X$.
In this case, the complex $\mathcal{K}$ together with the homeomorphism is called a \emph{triangulation}\index{triangulation}%
\footnote{The term is a bit misleading, as a triangulation may contain simplices of dimension larger than $2$.}
of~$X$.
The images of simplices of $\mathcal{K}$ in $X$ will be called the simplices of the triangulation\footnote{Since the images of straight lines under a homeomorphism are curves, the simplices in a topological polytope may not ``look like'' Euclidean simplices.}.

\begin{thm}{Definition}\label{def:poly}
A length space $P$ is called a \emph{polyhedral space}\index{polyhedral space}
if it admits a triangulation such that each simplex in $P$ is isometric to a simplex in Euclidean space.
\end{thm}

When we refer to a triangulation of a polyhedral space in the future, we will always mean a triangulation as in the above definition.

To construct an example of a polyhedral space,
one may take the underlying set of any simplicial complex in $\RR^n$
and equip it with the induced length metric.%
\footnote{In fact, any polyhedral space is isometric to one of these examples; see Exercise~\ref{ex:zalgalle+embedding}.}

Given a triangulation of a polyhedral space, 
we can consider the associated \emph{barycentric coordinates}\index{barycentric coordinates}.
If $\{v_1,\dots,v_n\}$ is the set of all vertices of the triangulation,
then any point $x$ is described uniquely by $n$ numbers $\lambda_1,\dots,\lambda_n$
such that $\lambda_i\ge 0$ for all $i$,
$\lambda_1+\dots+\lambda_n=1$,
and there is a simplex $\Delta$ in the triangulation such that
$\lambda_i\ne 0$ if and only if $v_i$ is a vertex of $\Delta$.
Indeed, let $\Delta$ be the minimal simplex which contains $x$.
Since $\Delta$ is isometric to a Euclidean simplex,
we can use barycentric coordinates described above and take $\lambda_i=0$
if $v_i$ is not a vertex of $\Delta$.

The \emph{dimension of a polyhedral space}\index{dimension of polyhedral space} is defined as the maximum dimension of the simplices in its triangulation.
(One can show that this value does not depend on the choice of triangulation.)

\begin{thm}{Exercise}
Prove that a convex polygon $A$ in $\RR^2$ with the subspace metric is a 2-dimensional polyhedral space.
Moreover, show that for any triangulation of $A$ as in the definition of polyhedral space, the simplices of $A$ are Euclidean simplices.
\end{thm}

Two other important examples of polyhedral spaces are
a convex polyhedron in $\RR^3$
and the boundary of a convex polyhedron equipped with its induced length metric.
In order to prove that these spaces are indeed polyhedral, 
one might triangulate these spaces by hand.
Along the same lines, one can also prove the following characterization of polyhedral spaces:

\begin{thm}{Theorem}
A length space $P$ is a polyhedral space if it can be covered by a finite number of subsets $M_1,\dots, M_n$
such that each $M_i$, as well as every intersection of a subcollection of $M_i$'s, is isometric to a convex polyhedron.
\end{thm}

\addtocontents{toc}{\protect\end{quote}}

\chapter{Zalgaller's folding theorem}
\addtocontents{toc}{\protect\begin{quote}}
\addtocontents{toc}{\sloppy Any polyhedral space admits a piecewise distance-preserving map into Euclidean space of the same dimension.}

Let $P$ be a polyhedral space.
A map $f\:P\to\RR^n$ is called 
\emph{piecewise distance-preserving}\index{piecewise distance-preserving map}
if there is a triangulation of $P$ 
such that for any simplex $\Delta$ in the triangulation, 
the restriction $f|_\Delta$ is distance-preserving.

\begin{thm}{Exercise}\label{ex:PDP=>cont}
Show that any piecewise distance-preserving map $f$ is continuous and length-preserving.
(It then follows from Exercise \ref{LP=>short} that $f$ is also distance non-expanding.)
\end{thm}

The following statement might look obvious, 
but try to prove it rigorously.

\begin{thm}{Exercise}\label{ex:n=<m}
Suppose that an $m$-dimensional polyhedral space admits a piecewise distance-preserving map into $\RR^n$.
Show that $n\ge m$.
\end{thm}

In fact, the converse of the statement in the exercise is true.
In other words, dimension is the only obstruction to the existence of a piecewise distance-preserving map into Euclidean space.
The following theorem asserts this for the case where $m=2$.

\begin{thm}{Zalgaller's theorem}\label{thm:zalgaller}
Any 2-dimensional polyhedral space admits a piecewise distance-preserving map into the Euclidean plane.
\end{thm}

Imagine that you have a paper model of a 2-dimensional polyhedral space $P$ in your hands, and you fold this model so that it lies flat on a table.%
\footnote{We recommend creating such a paper model, say the surface of a cube, and then trying to fold it on the table.}
This is an intuitive way to think of a piecewise distance-preserving map $f\: P \to \RR^2$.
To make it closer to the actual definition, one has to imagine that the layers of paper can go thru each other.
Zalgaller's theorem says that such a ``folding'' is always possible;
see also Exercise~\ref{pr:6-4-3-4}.

The following exercise shows that in this process, new folds may need to be introduced across the triangles of the given triangulation of $P$.

\begin{thm}{Exercise}\label{pdp-for-tetrahedron}
Let $\Delta$ be a non-degenerate 3-dimensional simplex in $\RR^3$
and let $\partial\Delta$ be its boundary, equipped with the induced length metric.
It is a polyhedral space glued from 4 triangles --- the faces of $\Delta$.

Show that $\partial\Delta$ does not admit a map to $\RR^2$ which is distance-preserving on each of the 4 faces of $\Delta$.

Describe explicitly a piecewise distance-preserving map
$$f\:\partial\Delta\z\to\RR^2$$
which is distance-preserving
on 2 out of the 4 faces of $\Delta$.
(You will have to subdivide the other 2 faces into smaller triangles.)
\end{thm}

Below we give two similar proofs of Zalgaller's theorem: 
the first with cheating and the second without.
In the first proof, we use the following claim without proof.

\begin{clm}{}\label{clm:acute-triangulation}
Any 2-dimensional polyhedral space admits an 
\emph{acute triangulation}\index{acute triangulation},
that is, a triangulation such that all of its triangles are acute.
\end{clm}

\begin{thm}{Exercise}\label{ex:acute-triangulation}
Show that any triangle admits an acute triangulation.
\end{thm}

A proof of \ref{clm:acute-triangulation} is given in \cite{saraf}.
The proof requires more than simply subdividing each triangle into acute triangles --- the subdivisions of two triangles that share an edge must be compatible.

\parit{Proof using \ref{clm:acute-triangulation}.}
Fix an acute triangulation $\mathcal{T}_0$ of $P$ provided by \ref{clm:acute-triangulation}.
Mark all its vertices in white
and denote them by $\{w_1,\dots,w_k\}$.

For each $w_i$, consider its \emph{Voronoi domain}\index{Voronoi domain} $V_i$,
which is the subset
$$V_i=\set{x\in P}{|x-w_i|\le|x-w_j|\ \text{for any}\ j}.$$
Denote by $S(w_i)$ the \emph{star}\index{star} of $w_i$,
which is the union of all simplices of the triangulation $\mathcal{T}_0$ which contain $w_i$.

Since an acute triangle contains its own circumcenter, it is impossible for the Voronoi domain of a vertex of a triangle to cross the opposite edge.
From this it follows that
$$V_i\subset S(w_i)$$ for all $i$.
In particular, for any point $x\in V_i$, there is a unique geodesic $[w_i,x]$,
which is a line segment in a single triangle or an edge of $\mathcal{T}_0$.

%Note that since all triangles are acute, we have
%$$V_i\subset S(w_i)$$ for all $i$.
%In particular, for any point $x\in V_i$, there is a unique geodesic $[w_i,x]$,
%which is a line segment in a triangle or an edge of $\mathcal{T}_0$.

\begin{wrapfigure}{r}{49mm}
\begin{lpic}[t(-0mm),b(-0mm),r(0mm),l(0mm)]{pics/zalgaller(1)}
\end{lpic}
\caption*{The Voronoi domains within one triangle.}
\end{wrapfigure}

Note that in each triangle of $\mathcal{T}_0$, we have one point where three Voronoi domains meet and three points on the sides of the triangle where pairs of Voronoi domains meet.
Let us bisect each edge of $\mathcal{T}_0$
and subdivide each triangle into $6$ triangles as shown in the picture (solid lines only).

In this way we obtain a new triangulation $\mathcal{T}_1$.
We mark all the new vertices of $\mathcal{T}_1$ in black.

Note that 
\begin{enumerate}
\item Each $V_i$ is a union of all triangles and edges of $\mathcal{T}_1$ which have $w_i$ as a vertex.
\item Each triangle in $\mathcal{T}_1$ has one white and two black vertices.
\item\label{prop:cong-pairs} The triangles in $\mathcal{T}_1$ come in pairs of congruent triangles; they share two black vertices and have different white vertices.
\end{enumerate}

Given a point $x\in P$,
set
$$\rho(x)=\min_i\{|w_i-x|\}.$$
Notice that if $x \in V_i$, then $\rho(x) = |w_i - x|$.
Given $x\in V_i$, 
we denote by $\theta_i(x)$ the minimum angle between $[w_i,x]$ and any edge of $\mathcal{T}_1$ coming from $w_i$.

By Property~\ref{prop:cong-pairs}, 
if $x\in V_i\cap V_j$ then $\theta_i(x)=\theta_j(x)$.
In other words, the function $\theta$ given by
$$\theta(x) = \theta_i(x), \qquad x \in V_i$$ is well-defined on the set $P\backslash\{w_1,\dots,w_n\}$.
Moreover, $\theta$ is a continuous function.

We now describe the map $f\:P\to\RR^2$ using polar coordinates on~$\RR^2$.
We define
$f(w_i)=0$ and 
$f(x)=(\rho(x),\theta(x))$ if $x \in P\backslash\{w_1,\z\dots,w_n\}$.

Subdividing each triangle by the angle bisector at the white vertex (see the dashed lines in the picture)
gives a new triangulation $\mathcal{T}_2$
which satisfies the conditions of the theorem for the constructed map $f$.
\qeds

Now we modify the above proof 
so it does not use Claim \ref{clm:acute-triangulation}.
The only property we really need for the triangulation is that $V_i\subset S(w_i)$.
This inclusion does not hold for a general triangulation.
A simple example of this is obtained by gluing an equilateral triangle to the longer side of an obtuse triangle.
On the other hand, 
by increasing the number of Voronoi domains, we can arrange that this inclusion holds without making the triangulation acute.

\parit{Proof without using \ref{clm:acute-triangulation}.}
Fix a triangulation $\mathcal{T}_0$ of $P$ that is not necessarily acute.
We will construct new triangulations $\mathcal{T}_1$ and $\mathcal{T}_2$ of $P$
by subdividing the triangles of $\mathcal{T}_0$,
and we will define a map $f\:P\to\RR^2$ which is distance-preserving on each triangle of $\mathcal{T}_2$.
The vertices of $\mathcal{T}_1$ will be colored either white or black 
in such a way that each triangle of $\mathcal{T}_1$ 
will have two black vertices and one white vertex.
%(In addition to the sides and vertices of the triangles, $\mathcal{T}_1$  might have white-black edges.) 

\begin{wrapfigure}{r}{56mm}
\begin{lpic}[t(-4mm),b(-0mm),r(0mm),l(0mm)]{pics/voronoi(1)}
\end{lpic}
\caption*{A triangle $\Delta$ of $\mathcal{T}_0$
with marked white points and the intersections of their Voronoi domains with $\Delta$.}
\end{wrapfigure}

We shall first describe the set of white vertices.

Fix a small number $\eps>0$.
We mark in white all of the vertices of $\mathcal{T}_0$,
as well as the points on the sides of triangles of $\mathcal{T}_0$ 
with the property that the distance to the closest vertex of the edge is an integer multiple of $\eps$.
In this way we mark a finite number of points white.
Label the white points by $w_1,\dots,w_k$.

As in the previous proof, let $V_i$ be the Voronoi domain of $w_i$, so that
\[V_i=\set{x\in P}{|x-w_i|\le|x-w_j|\ \text{for any}\ j}.\]
We also let $S(w_i)$ be the star of $w_i$ in $\mathcal{T}_0$, which is the union of all simplices of $\mathcal{T}_0$ which contain $w_i$.
By the following exercise, we can assume that $V_i\subset S(w_i)$ for each $i$, by taking a suitably small $\eps$.

\begin{thm}{Exercise}\label{ex:voronoi-in-star}
Let $\ell$ be the minimal length of the edges in the triangulation,
and let $\alpha$ be the minimal angle of all the triangles in $\mathcal{T}_0$.
Show that if $\eps<\tfrac{\ell\cdot \alpha}{100}$, 
then $V_i\subset S(w_i)$ for each $i$.
\end{thm}

Fix a triangle $\Delta$ of $\mathcal{T}_0$.
Note that for any $w_i\in \Delta$, the intersection $V_i\cap \Delta$ is a convex polygon.
This follows because for any two points $w_i,w_j\in\Delta$,
the inequality 
$$|x-w_i|\le |x-w_j|$$ 
describes the set of all points $x \in \Delta$
which lie on one side of the bisecting perpendicular to $w_i$ and $w_j$.
Let us color the vertices of all of the polygons $V_i\cap \Delta$ in black, if they are not already white.

If $\mathcal{T}_0$ contains an edge $E$ which is not
a side of a triangle, then color the midpoint of $E$ in black.

We'll now describe the triangulation $\mathcal{T}_1$.
The vertices of $\mathcal{T}_1$ are the black and white vertices.
A white point $w_i$ is connected by an edge to each black point in $V_i$.
The black vertices $b$ and $b'$ are connected by an edge if they form a side of some $V_i\cap \Delta$
(for some $\Delta$ and $w_i\in\Delta$).
In this case $w_i$, $b$ and $b'$ also form a triangle of $\mathcal{T}_1$.
Notice that each black-black edge is a side of two congruent triangles of $\mathcal{T}_1$ with different white vertices.

The remaining part of the proof is the same  as before.
We define $$\rho(x)=\min_i\{|w_i-x|_P\}$$
and $\theta(x)$ for $x\in V_i$ as the minimal angle between $[w_i,x]$ and any edge
in $\mathcal{T}_1$ coming from $w_i$.
Then define the map $f\:P\to \RR^2$ so that 
$f(w_i)=0$ for each $i$ 
and
$f(x)=(\rho(x),\theta(x))$ in polar coordinates.

Further subdividing each triangle of $\mathcal{T}_1$
into two along the angle bisector from the white vertex produces a new triangulation $\mathcal{T}_2$.
It is straightforward to see that the constructed map $f$ is distance-preserving on each triangle of $\mathcal{T}_2$.
\qeds

Use Zalgaller's theorem to show the following.

\begin{thm}{Advanced exercise}\label{ex:zalgalle+embedding}
Any 2-dimensional polyhedral space is isometric to the underlying set of a simplicial complex in $\RR^n$, equipped with its induced length metric.
\end{thm}

We end this section with an entertaining exercise.

\begin{thm}{Exercise}\label{ex:black-and-white}
Let $\mathcal{T}$ be a triangulation of a convex polygon $Q$ in $\RR^2$ such that each triangle is colored either black or white.
Show that the following two conditions are equivalent.
\begin{enumerate}[a)]
\item There is a piecewise distance-preserving map
$Q\to\RR^2$ for this triangulation which preserves the orientation%
\footnote{We say that the motion preserves/reverses the orientation 
if it is a composition of an even/odd number of reflections.}
of each white triangle
and reverses the orientation of each black triangle.

\item The sum of black angles around any vertex  of $\mathcal{T}$ 
which lies in the interior of $Q$
is either $0$, $\pi$ or $2\cdot\pi$.
\end{enumerate}

\end{thm}

\addtocontents{toc}{\protect\end{quote}}
\chapter{Brehm's extension theorem}
\addtocontents{toc}{\protect\begin{quote}}
\addtocontents{toc}{\sloppy Any distance non-expanding map from a finite set of points in Euclidean space to the same Euclidean space can be extended to the whole space as a piecewise distance-preserving map.}

\begin{thm}{Brehm's extension theorem}\label{thm:brehm}
Let $a_1,\dots,a_n$ and $b_1,\dots,b_n$ be two collections of points in $\RR^2$ such that 
$$|a_i-a_j|\ge |b_i-b_j|$$
for all $i$ and $j$,
and let $A$ be a convex polygon which contains $a_1,\z\dots,a_n$.
Then there is a piecewise distance-preserving map $f\:A\to \RR^2$
such that
$f(a_i)\z=b_i$ for all $i$.
\end{thm}

In other words, if $F=\{a_1,\dots,a_n\}$ is a finite subset of a convex polygon $A$, then any distance non-expanding map $\phi\: F \to \RR^2$ extends to a piecewise distance-preserving map $f\: A \to \RR^2$.

\parit{Proof.}
The proof is by induction on $n$.

The base case $n=1$ is trivial:
we can take
$$f(x) = x + (b_1 - a_1),$$ which is distance-preserving on $A$.

Applying the induction hypothesis to the last $n-1$ pairs of points,
we get a piecewise distance-preserving map $h\:A\z\to \RR^2$
such that $h(a_i)=b_i$ for all $i>1$.
We will use $h$ to construct the desired map $f\:A\to \RR^2$.

Consider the set 
$$\Omega=\set{x\in A}{|a_1-x|<|b_1-h(x)|}.$$
We can assume that $a_1 \in \Omega$;
otherwise, $h(a_1)= b_1$ and we can take $f = h$.
We make the following claim.

%First note the following.

\begin{clm}{}\label{clm:star-shaped}
The set $\Omega$ is star-shaped with respect to $a_1$.
That is, if $x\in \Omega$, then the line segment $[a_1,x]$ lies in $\Omega$.
\end{clm}

Indeed, if $y\in [a_1,x]$, then
$$|a_1-y|+|y-x|=|a_1-x|.
$$
Since $x\in\Omega$, we have
$$|a_1-x| < |b_1-h(x)|.
$$
Since $h$ is distance non-expanding (see Exercise~\ref{ex:PDP=>cont}),
we have
$$|h(x)-h(y)|\le |x-y|.
$$
Combining the above with the triangle inequality, we see
\begin{align*}
|a_1-y| &= |a_1-x| - |x-y|
<
\\
&< |b_1-h(x)| - |h(x)-h(y)|
\le
\\
&\le
|b_1-h(y)|.
\end{align*}
This proves $y\in\Omega$, which establishes Claim~\ref{clm:star-shaped}.

\medskip

\begin{center}
\begin{lpic}[t(-0mm),b(-0mm),r(0mm),l(0mm)]{pics/brehm-new(1)}
\lbl{55,55;\Large{$A$}}
\lbl[lb]{45,35;$E_i$}
\lbl[]{38,33;$T_i$}
\lbl[tr]{29,28;$a_1$}
\lbl[]{15,28;\Large{$\Omega$}}
\lbl[b]{10,35,25;\Large{$\partial_A\Omega$}}
\lbl[b]{-1,15,90;\Large{$Z$}}
\lbl[b]{15,10;blind}
\lbl[t]{15,9;zone}
\end{lpic}
\end{center}

Recall that $\partial_A\Omega$ denotes the boundary of $\Omega$ considered as a subset of the space $A$.
This may be a different set from $\partial_{\RR^2}\Omega$.
Note that  
$$|a_1-x|=|b_1-h(x)|
\eqlbl{eq:|ax|=|ah(x)|}$$
for any $x\in\partial_A\Omega$.
To see this, consider a sequence of points in $\Omega$ that converges to $x$ and another sequence of points in $A\backslash \Omega$ that converges to $x$, and then use the fact that $h$ is continuous (Exercise~\ref{ex:PDP=>cont}).

Further, note the following.
\begin{clm}{}\label{clm:broken-line}
The boundary $\partial_A\Omega$ is the union of a finite collection of line segments
$E_1,\dots, E_k$ which intersect only at common endpoints.
Moreover, $h$ is distance-preserving on each of these segments.
\end{clm}

Indeed, fix a triangulation of $A$ so that $h$ is distance-preserving on each triangle.
Note that for any point $x\in \partial_A\Omega$, this triangulation has a triangle $\Delta\ni x$ such that $\Delta\cap\Omega\ne\emptyset$.
Fix such a triangle $\Delta$.
Since $h$ is distance-preserving on $\Delta$,
the restriction $h|_\Delta$ can be extended uniquely to an isometry $\iota\:\RR^2\z\to\RR^2$.

Set $b_1'=\iota^{-1}(b_1)$.
Note that 
\[|b_1'-x|=|b_1-h(x)|\] 
for any $x\in\Delta$, because $\iota$ is an isometry and $\iota|_\Delta = h|_\Delta$.

Observe that $a_1\ne b_1'$.
Assuming otherwise, we see $$|a_1-x|=|b_1'-x|=|b_1-h(x)|$$ for any $x\in\Delta$, which gives the contradiction $\Delta\cap\Omega=\emptyset$.

Denote by $\ell_\Delta$ the perpendicular bisector of $[a_1, b_1']$, which coincides with the set of all points equidistant from $a_1$ and $b_1'$.
By the definition of $\Omega$, for any $x\in \Delta$ we have that
$x\in\Omega$ if and only if $x$ and $a_1$ lie on the same side of $\ell_\Delta$.
Therefore $\partial_A\Omega$ is the union of the intersections $\Delta\cap\ell_\Delta$ for all $\Delta$ as above.
Hence \ref{clm:broken-line} follows, as there are only finitely many such $\Delta$.

For each edge $E_i$ in  $\partial_A\Omega$, consider the triangle $T_i$ with vertex $a_1$ and base $E_i$.
Condition \ref{eq:|ax|=|ah(x)|} implies that there is an isometry $\iota_i$ of $\RR^2$ such that $\iota_i(a_1)= b_1$
and  $\iota_i(x)=h(x)$ for any $x\in E_i$.

Let us define 
$f(x)=h(x)$ for any $x\notin\Omega$
and $f(x)=\iota_i(x)$ for any $x\in T_i$.
This defines $f$ on $A\backslash \Omega$ 
and on all line segments from $a_1$ to~$\partial_A\Omega$.

This completely defines $f$ on $A$ in the case where $\partial_{A}\Omega=\partial_{\RR^2}\Omega$.
If $Z=\partial_{\RR^2}\Omega\backslash\partial_{A}\Omega$ is nonempty, then the points between $a_1$ and the points in $Z$
form a ``blind zone'' --- this is the subset of $A$ where $f$ has yet to be defined.

Note that the closure of the blind zone is a union of a finite number 
of polygons $Q_1,\dots, Q_m$ which intersect only at the common vertex $a_1$.
Each $Q_i$ is bounded by a broken line in the closure of $Z$ 
and two line segments from $a_1$ to the ends of this broken line.

So far the distance non-expanding map $f$
is defined only on the two sides of each $Q_i$ coming from $a_1$, 
and by construction it is distance-preserving on each of these two sides.
From the exercise below, it follows that
one can extend $f$ to each $Q_i$ while keeping it piecewise distance-preserving.
\qeds

\begin{thm}{Exercise}\label{ex:triangle-reflect}
Let $Q=[a_1x_1\dots x_k]$ be a polygon and $b_1$, $y_1$, $y_k$ be points in the plane.
Assume that 
\begin{align*}
|b_1-y_1|&=|a_1-x_1|,&
|b_1-y_k|&=|a_1-x_k|,&
|y_1-y_k|&\le|x_1-x_k|.
\end{align*}
Then there is a piecewise distance-preserving map
$f\:Q\to \RR^2$ such that $f(x_1)=y_1$, $f(x_k)= y_k$ and $f(a_1)= b_1$.
\end{thm}

Let us finish this lecture with some additional exercises.

\begin{thm}{Exercise}\label{pr:perimeter}
Let $a_1,\dots,a_n$ and $b_1,\dots,b_n$ be two collections of points in $\RR^2$ such that 
$$|a_i-a_j|\ge |b_i-b_j|$$
for all $i$ and $j$.
Let $A=\Conv\{a_1,\dots,a_n\}$ 
and $B=\Conv\{b_1,\z\dots,b_n\}$ be their convex hulls.
Show that 
$$\per A\ge \per B,$$
where $\per A$ denotes the perimeter of $A$.

Is it true that
$$\area A\ge \area B?$$

\end{thm}

The following exercise is a 2-dimensional case of Alexander's theorem
\cite{alexander}.
It has quite a simple solution,
but it plays an important role in discrete geometry;
check, for example, the paper by Bezdek and Connelly \cite{bezdek-connelly}.

\begin{thm}{Advanced exercise}\label{pr:alexander}
Let $a_1,\dots,a_n$ and $b_1,\dots,b_n$ be two collections of points in $\RR^2$.
Let us consider $\RR^2$ as a coordinate plane $\RR^2 \times \{0\}$
in $\RR^4=\RR^2 \times \RR^2$.

Construct a collection of curves $\alpha_i\:[0,1]\to \RR^4$ 
such that
$\alpha_i(0)\z=a_i\z=(a_i,0)$, $\alpha_i(1)=b_i=(b_i,0)$ and
the function $\ell_{i,j}(t)\z=|\alpha_i(t)\z-\alpha_j(t)|$
is monotonic (i.e., increasing, decreasing or constant) for each $i$ and $j$.
\end{thm}

\begin{thm}{Exercise}\label{pr:brehm}
Use Brehm's extension theorem 
to prove Kirszbraun's theorem, stated below, in the special case where $Q$ is a finite set.
\end{thm}

\begin{thm}{Kirszbraun's theorem}
Let $Q\subset \RR^2$ be an arbitrary subset and $f\:Q\to\RR^2$ be a distance non-expanding map.
Then $f$ admits 
a distance non-expanding extension to all of $\RR^2$.
In other words, there is a distance non-expanding map $F\:\RR^2\to\RR^2$
such that the restriction $F|_{Q}$ coincides with $f$.
\end{thm}

\addtocontents{toc}{\protect\end{quote}}
\chapter{Akopyan's approximation theorem}
\addtocontents{toc}{\protect\begin{quote}}
\addtocontents{toc}{Any distance non-expanding map from a polyhedral space to Euclidean space of the same dimension can be approximated by a piecewise distance-preserving map.}

Let $P$ be a polyhedral space.
A map $h\:P\to\RR^n$ is called 
\emph{piecewise linear}\index{piecewise linear map}\label{page:piecewise linear map}
if there is a triangulation of $P$ such that 
the restriction of $h$ to any simplex $\Delta$ is a linear map.
This means that if $v_0,\dots,v_k$ are the vertices of $\Delta$,
then for any $x\in \Delta$ we have
\[h(x)
=
\lambda_0\cdot h(v_0)+\dots+\lambda_k\cdot h(v_k),\]
where $(\lambda_0,\dots,\lambda_k)$ are the barycentric coordinates of $x$.

\begin{thm}{Exercise}\label{ex:PDPisPL}
Show that if $P$ is a $2$-dimensional polyhedral space, then any piecewise distance-preserving map $f\: P \to \RR^2$ is piecewise linear.
\end{thm}

In general, piecewise linear maps may not be injective, and they may either expand or contract distances.
We shall be interested in approximating piecewise linear maps by piecewise distance-preserving maps.
Since all piecewise distance-preserving maps are distance non-expanding, it only makes sense to try this approximation for distance non-expanding maps.

\begin{thm}{Akopyan's theorem}\label{thm:approx}
Assume $P$ is a $2$-dimensional polyhedral space.
Then any distance non-expanding piecewise linear map
$h\:P\z\to\RR^2$ can be approximated 
by piecewise distance-preserving maps.

More precisely, given $\eps>0$, there is a piecewise distance-preserving map $f\:P\to\RR^2$
such that 
$$|f(x)-h(x)|<\eps$$
for all $x\in P$.

\end{thm}

This theorem implies the existence of many piecewise distance-preserving maps from $P$ into $\RR^2$.
In particular, it implies Zalgaller's theorem (\ref{thm:zalgaller}).
To see this, consider the constant map $h\:P\z\to\RR^2$;
i.e., the map which sends the whole space $P$ to a single point.
Since $h$ is piecewise linear, 
we can apply Akopyan's theorem to produce a piecewise distance-preserving map $f$ arbitrarily close to $h$.

As you will see below, 
the proof of Akopyan's theorem will not use Zalgaller's theorem.
We still consider the proof of Zalgaller's theorem to be important because it gives a very clear geometric description of the piecewise distance-preserving map.
In contrast, the maps produced by Akopyan's theorem will rely on the recursive construction of Brehm's theorem, which is harder to understand.

\begin{thm}{Exercise}\label{ex:akopyan-brehm}
Show that if $P$ is a convex polygon in $\RR^2$,
then the above theorem follows from Brehm's extension theorem (\ref{thm:brehm}).
\end{thm}

The main idea in the proof of Akopyan's theorem is to triangulate $P$ and use Brehm's extension theorem on each triangle, as in the previous exercise.
Unfortunately, it is not that simple.
The big technical issue that arises is that if two triangles share a common edge, then we need to ensure that the maps produced using Brehm's theorem agree on that common edge.

\begin{wrapfigure}{r}{43mm}
\begin{lpic}[t(-3mm),b(-2mm),r(0mm),l(6mm)]{pics/zig-zag(1)}
\lbl[r]{1,46;$h(z_0)$}
\lbl[r]{8,35;$h(z_1)$}
\lbl[r]{11,29;$h(x)$}
\lbl{18,18,-57;$\dots$}
\lbl[r]{28,1;$h(z_n)$}
\lbl[bl]{17,37;$w_n(x)$}
\end{lpic}
\end{wrapfigure}

To address this issue, we will use the following \emph{zigzag construction}\index{zigzag construction}.
It produces a piecewise distance-preserving map
which is close to
a given distance non-expanding linear map defined on a line segment.
For the construction, we fix a unit vector $e$ in $\RR^2$.
The choice of $e$ does not matter, but the same $e$ must be used uniformly in all zigzag constructions that follow.

\parit{Zigzag construction.} 
Let $E$ be a line segment and
$h\:E\to\RR^2$ be a distance non-expanding linear map.
Let $\ell = \length E$ and $\ell' = \length h(E)$.
Since $h$ is distance non-expanding, we have $\ell' \le \ell$.

Fix a positive integer $n$, and subdivide $E$ into $n$ equal intervals.
Denote by $z_0,\dots,z_n$ the endpoints of these intervals.

Note that the image $h(E)$ is either a line segment or a point.
In the first case, let $u$ be a unit normal vector to $h(E)$;
otherwise, let $u=e$.

Given $x\in E$, 
set 
\begin{align*}
s_n(x)&=\min_i\{|z_i-x|\},
\\
w_n(x)&=k\cdot s_n(x)\cdot u+h(x),
\end{align*}
where
$k=\sqrt{1-(\ell'/\ell)^2}$.
If we subdivide $E$ further by adding the midpoints between any two consecutive endpoints, then $w_n$ is distance-preserving on each of the resulting subintervals.
This shows that $w_n$ is piecewise distance-preserving.
Moreover
$$|w_n(x)-h(x)|\le\tfrac{\ell}{2\cdot n}$$
for any $x\in E$, because $k \le 1$ and $s_n(x) \le \tfrac{\ell}{2\cdot n}$.

The piecewise distance-preserving map $w_n$ is the result of the \emph{$n$-step zigzag construction} applied to $h$.

\medskip

Given a triangulation $\mathcal{T}$ of a polyhedral space $P$, let $\mathcal{T}^1$ denote the $1$-skeleton of $\mathcal{T}$.
This is the 1-dimensional subcomplex of $\mathcal{T}$
formed by all the vertices and edges in $\mathcal{T}$.
Notice that $\mathcal{T}^1$ is a $1$-dimensional polyhedral space when equipped with its induced length metric, which is different from the subspace metric it inherits from $P$.

The following proposition is the main technical step in the proof of Akopyan's theorem.

\begin{thm}{Proposition}\label{clm:main-step} 
Let $\mathcal{T}^1$ be the 1-skeleton of a triangulation of a $2$-dimensional polyhedral space $P$,
and let $h\:\mathcal{T}^1\to \RR^2$
be a piecewise linear map
such that 
$$|h(x)-h(y)|_{\RR^2}\le |x-y|_P$$
for any $x,y\in \mathcal{T}^1$.
Then for any $\eps > 0$, there is a piecewise distance-preserving map $w\:\mathcal{T}^1\to \RR^2$ such that
$$|w(x)-w(y)|_{\RR^2}\le |x-y|_P$$ for any $x,y \in \mathcal{T}^1$ and
$$|w(x) - h(x)| < \eps$$ for all $x \in \mathcal{T}^1$.
\end{thm}

\parit{Proof.}
First we prove the statement under the following additional assumption on $h$:

\begin{clm}{}\label{clm:delta-condition}
For some fixed $\delta>0$, we have
\[|h(x)-h(y)|_{\RR^2}\le(1-\delta)\cdot|x-y|_P\]
for any $x,y\in \mathcal{T}^1$
and 
$$h(v)=h(x)$$ 
for any vertex $v$ of $\mathcal{T}^1$
and any point $x\in\mathcal{T}^1$ such that  $|v-x|_P\le\delta$.
\end{clm}

Let $\mathcal{S}$ 
denote the subdivision of $\mathcal{T}^1$
such that $h$ is linear on each edge of $\mathcal{S}$.
Subdividing $\mathcal{S}$ further if necessary, we may assume without loss of generality that 
each edge of $\mathcal{S}$ 
which comes from a vertex of $\mathcal{T}^1$ 
has length $\delta$.
(To perform this subdivision, 
we have to assume that $\delta$ in \ref{clm:delta-condition} is sufficiently small.)

Denote by $\ell$ the maximal length of the edges in $\mathcal{T}^1$.
Let us apply the $n$-step zigzag construction to each edge of $\mathcal{S}$.
Since the maps from the zigzag construction agree at the common vertices of different edges, we obtain a piecewise distance-preserving map $w_n\:\mathcal{T}^1\to\RR^2$ such that
$$|w_n(x) - h(x)| \le \tfrac{\ell}{2 \cdot n}
\eqlbl{eq:wn=f}$$ 
for all $x \in \mathcal{T}^1$.

We shall show that the inequality
$$|w_n(x) - w_n(y)|_{\RR^2} \le |x - y|_P$$ holds for all $x,y \in \mathcal{T}^1$, provided $n$ is sufficiently large.
Notice that
\begin{clm}{}\label{eq:w-on-edge}
$|w_n(x)-w_n(y)|\le |x-y|_P$
if $x$ and $y$ lie on the same edge.
\end{clm}
\noi
Indeed, if $x$ and $y$ lie on the same edge $E$ of $\mathcal{T}^1$,
then $|x-y|_P = |x-y|_{E}$, and the map $w_n$ is distance non-expanding on $E$.

From \ref{clm:delta-condition} and \ref{eq:wn=f},  we see that 
\begin{align*}
&|w_n(x) - w_n(y)|_{\RR^2}\le\\ 
&\qquad \le |w_n(x) - h(x)|_{\RR^2} ~+~ |h(x) - h(y)|_{\RR^2} ~+~ |h(y) - w_n(y)|_{\RR^2}\le\\
&\qquad  \le |x-y|_P+\left(\tfrac{\ell}{n} - \delta\cdot|x-y|_P\right)
\end{align*}
for any $x$ and $y$ in $\mathcal{T}^1$.

Now suppose $|w_n(x)-w_n(y)|_{\RR^2}> |x-y|_P$ for some $x,y \in \mathcal{T}^1$.
Then from above, we have $|x - y|_P < \tfrac{\ell}{n\cdot \delta}$, which shows that
\[|x-y|_P<\tfrac{C}{n}\eqlbl{eq:x-near-y}\]
for a constant $C$ which does not depend on $x$ or $y$.
Thus, \ref{eq:w-on-edge} and \ref{eq:x-near-y} imply the following.

\begin{clm}{}\label{clm:near-vertex}
For sufficiently large%
\footnote{$n$ does not depend on $x$ or $y$.
To ensure $x$ and $y$ do not lie on disjoint edges, $n$ must be large so that $C/n$ is less than the minimal distance between any two disjoint edges.
To ensure that both $x$ and $y$ are within $\delta$ of $v$, we must choose $n$ large in a way which will depend on the minimal angle in any triangle.
To deal with both issues, we are using the fact that $\mathcal{T}$ has a finite number of triangles.} $n$,
if $|w_n(x)-w_n(y)|> |x-y|_P$ then both $x$ and $y$ 
lie on different edges which come from one vertex, say $v$ of $\mathcal{T}^1$,
and 
\[|x-v|_{P},\  |y-v|_{P}\le\delta.\] 
\end{clm} 
Let $x$, $y$ and $v$ be as in \ref{clm:near-vertex}.
Take the point $x'$ on the same edge as $y$
such that $|v-x'|_P=|v-x|_P$.
It follows from the construction of $w_n$ that
$w_n(x')=w_n(x)$.
(Notice that by \ref{clm:delta-condition}, $w_n$ is produced by the zigzag construction in the case where the image of $h$ is a point.)  Therefore
\begin{align*}
|w_n(x)-w_n(y)|_{\RR^2}&=|w_n(x')-w_n(y)|_{\RR^2}\le
\\
&\le |x'-y|_P=
\\
&=\bigl||x-v|_P-|y-v|_P\bigr|\le
\\
&\le|x-y|_P.
\end{align*}

Thus we have shown that if $n$ is sufficiently large, then the inequality 
$$|w_n(x)-w_n(y)|_{\RR^2}\le |x-y|_P$$
holds for any pair $x,y\in\mathcal{T}^1$.
Let $w = w_n$ for such an $n$ which additionally satisfies $\tfrac{\ell}{2\cdot n} < \eps$.
Then from \ref{eq:wn=f}, it follows that $$|w(x) - h(x)| < \eps$$ for all $x \in \mathcal{T}^1$.

\begin{wrapfigure}{r}{45mm}
\begin{lpic}[t(-8mm),b(0mm),r(0mm),l(0mm)]{pics/q-graph(1)}
%\lbl[t]{20,0;$x$}
%\lbl[r]{1,18;$q_n(x)$}
\lbl[b]{6.5,8;$\delta$}
\lbl[b]{39,8;$\delta$}
\end{lpic}
\caption*{The graph of $q_\delta$ on one edge.}
\bigskip
\end{wrapfigure}

Thus we have proved the proposition under the assumption \ref{clm:delta-condition}.
It remains to be shown that the general case can be reduced to the case where
\ref{clm:delta-condition} holds.
We shall achieve this by approximating $h$ by a map that satisfies~\ref{clm:delta-condition}.

For small $\delta > 0$ (say less than half the smallest edge length), consider the map 
$$q_\delta\:\mathcal{T}^1\to\mathcal{T}^1$$ 
which smashes the $\delta$-neighborhood of each vertex of $\mathcal{T}^1$ to the vertex 
and linearly stretches the remaining part of the edge, as in the figure.

Let $L_\delta$ be the optimal Lipschitz constant of $q_\delta$;
i.e., the minimal number such that 
\[|q_\delta(x)-q_\delta(y)|_P\le L_\delta\cdot|x-y|_P\] for all $x,y \in \mathcal{T}^1$.
Notice that $L_\delta\z\to 1$ as $\delta\to 0^+$.
Then the map 
\[h_\delta
\df
\tfrac{1-\delta}{L_\delta}\cdot (h\circ q_\delta)\]
is piecewise linear and satisfies condition \ref{clm:delta-condition}.
Moreover, we can choose $\delta$ sufficiently small so that
\[|h_\delta(x) - h(x)| < \tfrac{\eps}{2}\]
for all $x \in \mathcal{T}^1$.

By the previous part of the proof, there is a piecewise distance-preserving map $w\: \mathcal{T}^1 \to \RR^2$ such that
$$|w(x) - h_\delta(x)| < \tfrac{\eps}{2}, \qquad |w(x) - w(y)|_{\RR^2} \le |x - y|_P$$
for all $x, y \in \mathcal{T}^1$.
By the triangle inequality,
$$|w(x) - h(x)| < \eps$$ for all $x \in \mathcal{T}^1$.
\qeds

\parit{Proof of \ref{thm:approx}.}
Fix a fine
triangulation $\mathcal{T}$
of $P$,
one for which the diameter of each triangle is smaller than $\tfrac\eps{3}$.
Let $\mathcal{T}^1$ denote the $1$-skeleton of~$\mathcal{T}$.
By Proposition~\ref{clm:main-step}, 
there is a piecewise distance-preserving map $w\:\mathcal{T}^1\to\RR^2$ such that
$$|w(x)-h(x)|_{\RR^2} < \tfrac\eps{3}$$
for any $x\in \mathcal{T}^1$ and
$$|w(x)-w(y)|_{\RR^2}\le|x-y|_P$$
for any $x$ and $y\in \mathcal{T}^1$.

We shall use Brehm's extension theorem (\ref{thm:brehm}) to extend $w$ to a piecewise distance-preserving map on $P$.
To do this, let $\mathcal{S}$ be a subdivision of $\mathcal{T}^1$
so that $w$ is distance-preserving on each  edge of $\mathcal{S}$.
Fix a triangle $\Delta$ of $\mathcal{T}$.
Let $a_1,\dots, a_n$ be the vertices of $\mathcal{S}$ on the boundary of $\Delta$, 
and let $b_i=w(a_i)$ for each $i$.
By applying Brehm's theorem, we obtain a piecewise distance-preserving map $f_\Delta\: \Delta \to \RR^2$.

Since $w$ is distance-preserving on each edge of $\mathcal{S}$,
the maps $f_\Delta$ and $w$ coincide on the boundary of $\Delta$.
In particular, if $\Delta$ and $\Delta'$ share a common edge, then $f_\Delta$ and $f_{\Delta'}$ agree on that common edge.
Therefore the collection of maps $\{f_\Delta\}$ determines a single piecewise distance-preserving map $f\:P\to\RR^2$.

{\sloppy

We'll show that $f$ satisfies the conclusion of the theorem.
Let $x \in P$ be arbitrary and let $y$ be a point on the edge of a triangle in $\mathcal{T}$ that contains $x$.
Then $|x - y| < \tfrac{\eps}{3}$ by our choice of $\mathcal{T}$.
We see
\begin{align*}
|f(x) - h(x)| &\le |f(x) - w(y)| + |w(y) - h(y)| + |h(y) - h(x)|=\\
&= |f(x) - f(y)| + |w(y) - h(y)| + |h(y) - h(x)|\le\\
&\le 2\cdot |x - y| + |w(y) - h(y)|<\\
& < \eps,
\end{align*} 
because $w(y) = f(y)$ and both maps $f$ and $h$ are distance non-expanding.
\qeds

}

We close this section with a counterexample explaining one way in which we cannot improve Akopyan's theorem.
One might expect that a stronger statement holds, 
namely that the map $f$ in Akopyan's theorem 
can be constructed so that it coincides with $h$ 
on a given finite set of points.
The following exercise shows that this cannot be done in general.

\begin{thm}{Exercise}\label{ex:tripod}
Consider the following 5 points in $\RR^3$:
\begin{align*}
o=(0,0,0),
\,
p=(0,0,1),
\,
a=(2,0,0),
\,
b=(-1,2,0),
\,
c=(-1,-2,0)
\end{align*}

\begin{wrapfigure}{r}{37mm}
\begin{lpic}[t(-0mm),b(-0mm),r(0mm),l(0mm)]{pics/3-pod(1)}
\lbl[t]{15,10;$o$}
\lbl[b]{15,23.5;$p$}
\lbl[t]{35,10;$a$}
\lbl[rb]{7,23;$b$}
\lbl[rb]{0,2;$c$}
\end{lpic}
\end{wrapfigure}

Let $P$ be the ``tripod'' which is the polyhedral space consisting of the three triangles $\triangle opa$, $\triangle opb$ and $\triangle opc$ in $\RR^3$
and equipped with the induced length metric.

Note that the restriction of the coordinate projection $\pi(x,y,z) = (x,y,0)$ to $P$ is distance non-expanding and piecewise linear.
We have that 
\begin{align*}
\pi(o)&=\pi(p)=o,&\pi(a)&=a,&\pi(b)&=b,&\pi(c)&=c.
\end{align*}

Show that there is no piecewise distance-preserving map $f\:P\z\to \RR^2=\RR^2\z\times\{0\}$
such that $f(a)=a$, $f(b)=b$ and $f(c)=c$.
\end{thm}

\addtocontents{toc}{\protect\end{quote}}
\chapter{Gromov's rumpling theorem}\label{sec:S^2->R^2}
\addtocontents{toc}{\protect\begin{quote}}
\addtocontents{toc}{There is a length-preserving map from the sphere to the plane.}

Recall that $\SS^2$ denotes the unit sphere in $\RR^3$, which we equip with its induced length metric.
Here is our main theorem.

\begin{thm}{Theorem}\label{thm:S2->R2}
There is a length-preserving map $f\:\SS^2\to\RR^2$.
\end{thm}

Such a map $f$ has to crease on an everywhere dense set in $\SS^2$.
More precisely,  the restriction of $f$ to any open subset of $\SS^2$ cannot be injective.%
\footnote{To prove this, one can show that if $f$ is injective and length-preserving on an open set $U \subset \SS^2$, then $f$ maps (sufficiently short) geodesics to straight lines (this requires the Domain Invariance Theorem;
see, for example, Section 2.9 in \cite{alexandrov}).
It follows that the restriction of $f$ to $U$ is locally distance-preserving, which is impossible.
}

In the proof of the theorem we will use the following exercise.

\begin{thm}{Exercise}\label{problem2}
Let $K$ be a convex polyhedron in $\RR^3$.
Given a point $x$ in $\RR^3$, show that there is a unique point $\bar x\in K$ which minimizes the distance $|x-\bar x|$.
Moreover, show that the projection map $$\phi\: \RR^3 \to K, \qquad \phi(x) = \bar x$$ is distance non-expanding.
\end{thm}

\parit{Proof of Theorem~\ref{thm:S2->R2}.}
Consider a nested
sequence $K_0\subset K_1\subset \dots$ of convex polyhedra in $\RR^3$ whose union is the open unit ball.
Let $P_n = \partial K_n$ denote the surface of $K_n$, equipped with the induced length metric.
Note that $P_n$ is a $2$-dimensional polyhedral space for each $n$.

Let $\phi_n$ denote the distance non-expanding projection onto $K_n$
from Exercise~\ref{problem2}.
Since $K_{n} \subset K_{n+1}$, it follows that $\phi_n(P_{n+1}) = P_{n}$.
Note that one can triangulate $P_n$ and $P_{n+1}$ in such a way that the restriction of $\phi_n$ to any simplex of $P_{n+1}$ is an orthogonal projection to some simplex of $P_{n}$.
In particular, the restriction of $\phi_n$ to $P_{n+1}$ is piecewise linear%
\footnote{See the definition on page \pageref{page:piecewise linear map}.}
and distance non-expanding with respect to the length metrics on $P_{n+1}$ and $P_n$.

We claim that for any point $x\in \SS^2$, there is a unique sequence of points $x_n\in P_n$ such that $x_n\to x$ as $n\to\infty$ and $\phi_n(x_{n+1})=x_n$ for all $n$.
The uniqueness follows since the maps $\phi_n$ are distance non-expanding.
To show existence, fix any sequence $z_n\in P_n$ such that $z_n\to x$.
Consider the double sequence $y_{n,m}\in P_n$, defined for $n\le m$, such that $y_{n,n}=z_n$ and
$y_{n,m} = \phi_{n}(y_{n+1,m})$ if $0 \le n < m$.
Then set 
$$x_n=\lim_{m\to\infty} y_{n,m}.$$

\begin{thm}{Exercise}\label{ex:limit-above}
Show that the limit above exists
and $\phi_n(x_{n+1})\z=x_n$ for any $n$.
Then show that $x_n \to x$ as $n \to \infty$.
\end{thm}

Let $x_n\to x\in \SS^2$ be the sequence as above.
Define $\psi_n\:\SS^2\to P_n$ by $\psi_n(x)=x_n$.
We have that $\psi_n$ is distance non-expanding,
$\psi_n\z=\phi_n\circ\psi_{n+1}$ for all $n$,
and for any $p, q\in\SS^2$,
$$|p_n-q_n|_{P_n}\to |p-q|_{\SS^2}\ \text{as}\  n\to\infty,\eqlbl{eq:pq-to-pq}$$
where $p_n=\psi_n(p)$ and $q_n=\psi_n(q)$.

The desired length-preserving map $f\: \SS^2 \to \RR^2$ will be, in some sense, a ``limit'' of a sequence of piecewise distance-preserving maps $f_n\: P_n \to \RR^2$.
The maps will be constructed recursively to satisfy
$$|f_{n+1}(x) - f_n(\phi_n(x))| < \eps_n$$ for a carefully chosen sequence $(\eps_n)$ of positive numbers that decays rapidly to $0$.

\parit{Recursive construction of $f_n\:P_n\to\RR^2$ and $\eps_n$.}
Assume we have a piecewise distance-preserving map $f_n\: P_n\z\to \RR^2$ and a given $\eps_n$.
The composition $f_n\circ\phi_n\: P_{n+1}\z\to\RR^2$ is piecewise linear and distance non-expanding.
So we can apply Akopyan's theorem (see \ref{thm:approx}) to construct a piecewise distance-preserving map $f_{n+1}\:P_{n+1}\to\RR^2$ which is $\eps_n$-close to $f_n\circ\phi_n$.

Let $M(n+1)$ denote the number of triangles in a triangulation of $P_{n+1}$ such that $f_{n+1}$ is distance-preserving on each triangle.
Set 
$$\eps_{n+1}=\frac{\eps_n}{2\cdot M(n+1)}.\eqlbl{eq:eps-n}$$
In this way, we recursively define $f_n$ and $\eps_n$.
It goes as follows: 
\begin{enumerate}
\item Choose an arbitrary $\eps_0>0$ and a piecewise distance-preserving map $f_0\:P_0\to\RR^2$,
say the one provided by Zalgaller's folding theorem (\ref{thm:zalgaller}).
\item Use $\phi_0$, $f_0$ and $\eps_0$ to construct $f_1$.
\item Use $f_1$ to construct $\eps_1$.
\item Use $\phi_1$, $f_1$ and $\eps_1$ to construct $f_2$.
\item Use $f_2$ to construct $\eps_2$.
\item and so on.\footnote{The procedure is similar to walking on stairs: 
you take a right step, 
which makes it possible to take a left step, 
which in turn makes it possible to take a right step, and so on.}
\end{enumerate}

\medskip

It remains to prove the following claim:

\begin{clm}{}\label{clm:length-preserving}
The sequence of maps $f_n\circ\psi_n\:\SS^2\to \RR^2$ converges to a length-preserving map $f\:\SS^2\to\RR^2$.
\end{clm}

Since $\eps_{n}$ decays faster than $\tfrac{\eps_0}{2^n}$, 
the sequence $(f_n\circ\psi_n)(x) \in \RR^2$ is Cauchy, hence convergent,
for any fixed $x$.
We define $f\: \SS^2 \to \RR^2$ by
$$f(x) = \lim_{n \to \infty} (f_n\circ\psi_n)(x).$$
By the recursive construction of $f_n$, we have that
$$|(f_n\circ\psi_n)(x)-f(x)|<\eps_n$$
for any $x\in \SS^2$ and any $n$.
Since each $f_n\circ\psi_n$ is distance non-expanding, $f$ is distance non-expanding as well.

It only remains to show that the constructed map $f\: \SS^2 \to \RR^2$ is length-preserving.
By Exercise~\ref{LP=>short}(\ref{LP=>short:b}),
it suffices to show that
$$\length (f\circ\alpha) \ge |p-q|_{\SS^2}\eqlbl{length>=dist}$$
for any curve $\alpha$ between two points $p,q \in \SS^2$.
For the remainder of the proof, we will need the following definition;
it should be considered as an analog of the length of a curve.

\begin{thm}{Definition}
Let $X$ be a metric space and $\alpha\:[a,b]\to X$
be a curve.
Set 
$$\ell_k(\alpha)
\df
\sup\set{\sum_{i=1}^k |\alpha(t_i)-\alpha(t_{i-1})|_X}{a=t_0<t_1<\dots<t_k=b}.$$

\end{thm}

Note that given a curve $\alpha\:[a,b]\to X$,
we have
\begin{align*}
\ell_1(\alpha)&\le\ell_2(\alpha)\le \ell_3(\alpha)\le\dots,
\\
\ell_k(\alpha)&\to\length\alpha\ \ \text{as}\ \ k\to\infty,
\\
\ell_k(\alpha)&\le\length\alpha \ \text{for any}\ \ k.
\end{align*}
Moreover, if 
\[\ell_k(\alpha)=\length\alpha\] 
then $\alpha$ 
is a chain made from at most $k$ geodesic segments.

The following exercise states that if two curves $\alpha$ and $\beta$ are sufficiently close, then $\ell_k(\alpha)\approx\ell_k(\beta)$.
Note that according to the remark after \ref{thm:length-semicont},
the value $|\length\alpha-\length\beta|$ might be large in this case.

\begin{thm}{Exercise}\label{ex:klength-approx}
Suppose that $\alpha, \beta: \II \to X$ are two curves which are close in the sense that 
\[|\alpha(t) - \beta(t)|_X < \eps\] for all $t \in \II$.
Show that
\[ |\ell_k(\alpha) - \ell_k(\beta)| 
\le 2\cdot k\cdot \eps.\]

\end{thm}

Now we come back to the proof of \ref{length>=dist}.
Set $p_n=\psi_n(p)$ and $q_n\z=\psi_n(q)$.
Let $\beta$ be an arbitrary curve from $p_n$ to $q_n$ in $P_n$.
By Exercise~\ref{problem2}, one can find a shorter curve $\gamma$ from $p_n$ to $q_n$ whose image in any triangle of the triangulation of $P_n$ is a line segment, and moreover the endpoints of these line segments lie on $\beta$.
It follows that $f_n\circ\gamma$ is a broken line in $\RR^2$
with at most $M(n)$ edges, whose vertices we
denote, in order, by $$f_n(p_n) = z_0, z_1, \dots, z_k = f_n(q_n).$$
Note that  $k\le M(n)$
and 
each $z_i$ lies on the curve $f_n \circ \beta$.
%$z_i= (f_n(\beta(t_i))$ for some $a=t_0<t_1<\dots<t_k=b$.
Therefore
$$
\begin{aligned}
|p_n-q_n|_{P_n}
&\le \length\gamma=
\\
&= \ell_{M(n)}(f_n\circ\gamma)=
\\
&=|z_0-z_1|+\dots+|z_{k-1}-z_{k}|\le
\\
&\le\ell_{M(n)}(f_n\circ\beta);
\end{aligned}
\eqlbl{eq:pn-qn}
$$

Fix a curve $\alpha$ from $p$ to $q$ in $\SS^2$.
By Exercise \ref{ex:klength-approx} and \ref{eq:eps-n},
for all $n$ we have
$$|\ell_{M(n)}(f\circ\alpha)-\ell_{M(n)}(f_n\circ\psi_n\circ\alpha)| 
\le 2\cdot M(n)\cdot \eps_n = \eps_{n-1}.
\eqlbl{eq:inq-ell-k}$$
Given $\eps > 0$, 
we can choose $n$ sufficiently large
so that $\eps_{n-1} \le \tfrac{\eps}{2}$ and
$$|p-q|_{\SS^2} - |p_n - q_n|_{P_n} \le \tfrac{\eps}{2},$$ 
which can be arranged by \ref{eq:pq-to-pq}.
Applying \ref{eq:pn-qn} for $\beta=\psi_n \circ \alpha$ 
and \ref{eq:inq-ell-k}, we see
\begin{align*}
\length (f\circ \alpha) &\ge \ell_{M(n)}(f\circ \alpha)\ge\\
&\ge \ell_{M(n)}(f_n \circ \psi_n \circ \alpha) - \eps_{n-1}\ge\\
&\ge |p_n - q_n|_{P_n} - \eps_{n-1}\ge\\
&\ge |p-q|_{\SS^2} - \tfrac{\eps}{2} - \eps_{n-1}\ge\\
&\ge |p-q|_{\SS^2} - \eps.
\intertext{Since $\eps > 0 $ was arbitrary, }
\length(f\circ \alpha) &\ge |p-q|_{\SS^2}.
\end{align*}
Hence \ref{length>=dist} follows.
\qeds

\addtocontents{toc}{\protect\end{quote}}
\chapter[Arnold's problem on paper folding]{Arnold's problem\\ on paper folding}\label{sec:arnold}
\addtocontents{toc}{\protect\begin{quote}}
\addtocontents{toc}{Is it possible to fold a square in the plane so that the resulting figure will have a longer perimeter?}

This lecture is meant to be entertaining.
Here we will discuss the following problem posed by V.~Arnold in 1956 \cite[Problem 1956-1]{arnold}.

\begin{thm}{Problem}
Is it possible to fold a square in the plane so that the resulting figure will have a longer perimeter?
\end{thm}

\begin{wrapfigure}{r}{55mm}
\noi\begin{lpic}[t(-0mm),b(-5mm),r(0mm),l(0mm)]{pics/skladka(0.6)}
\lbl{25,28;$M$}
\lbl{75,28;$M'$}
\lbl[r]{23,55;$q$}
\lbl[r]{87,55;$q$}
\end{lpic}
\end{wrapfigure}

The answer to this problem depends on the meaning of the word ``fold''.

For example, one can consider a sequence of \emph{foldings} in which all layers are folded simultaneously along a line.
By the following exercise, the perimeter can never increase under a folding of this type.

\begin{thm}{Exercise}
Show that each fold described above indeed decreases the perimeter.
(Note that in general, the intersection of the line $q$ with the polygon $M$ in the picture may 
be a union of line segments.) 
\end{thm}

Using only the \emph{foldings} described above makes it impossible to unfold a layer which lies on top of another layer, as shown in the following picture.
\begin{center}
\begin{lpic}[t(0mm),b(0mm),r(0mm),l(0mm)]{pics/otgib(0.7)}
\end{lpic}
\end{center}
Note that this \emph{unfolding} increases the perimeter, altho not beyond the perimeter of the original square.
It is still unknown whether it is possible to increase the perimeter by a sequence of such ``folds'' and ``unfolds''.

\parbf{Japanese crane.}
Now let us consider a more general definition of folding.
Imagine that we mark in advance the lines of folding and start to fold the paper in such a way that each domain between folds remains flat all the time.

If you understand ``folding'' this way, then the answer to the problem is ``yes''.
In some sense, this problem was solved by origami practitioners well before it was even posed.

The possibility of increasing the perimeter slightly can be seen in the base for the crane.
This was known to origami masters for centuries%
\footnote{It appears in the oldest known book on origami, ``Senbazuru Orikata,'' dated 1797; but for sure it was known much earlier.},
but mathematicians learned this answer only in 1998%
\footnote{Here is \textattachfile[color=0 0 1,mimetype=text/html]{pics/napkin.html}{the html-file} which tells how it happened.}.

\begin{figure}[h]
\ \ \ \ \ \ \ \ 
\includegraphics[scale=0.34]{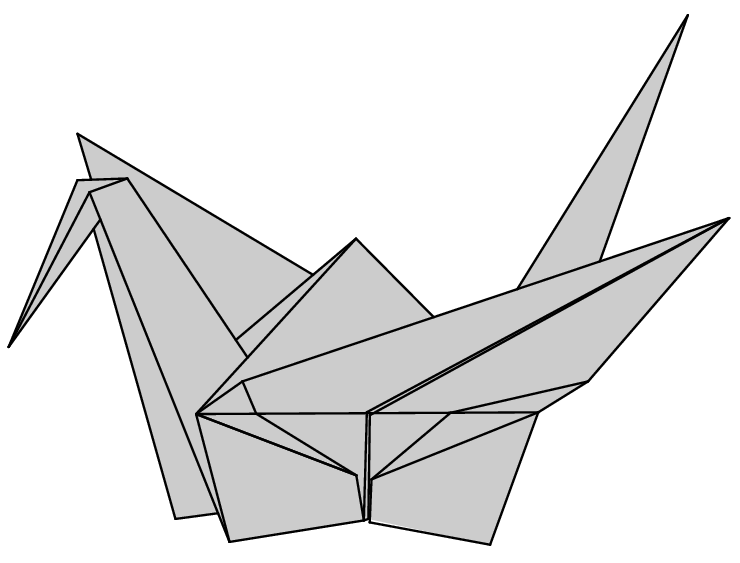}
\hfill
\includegraphics[scale=0.18]{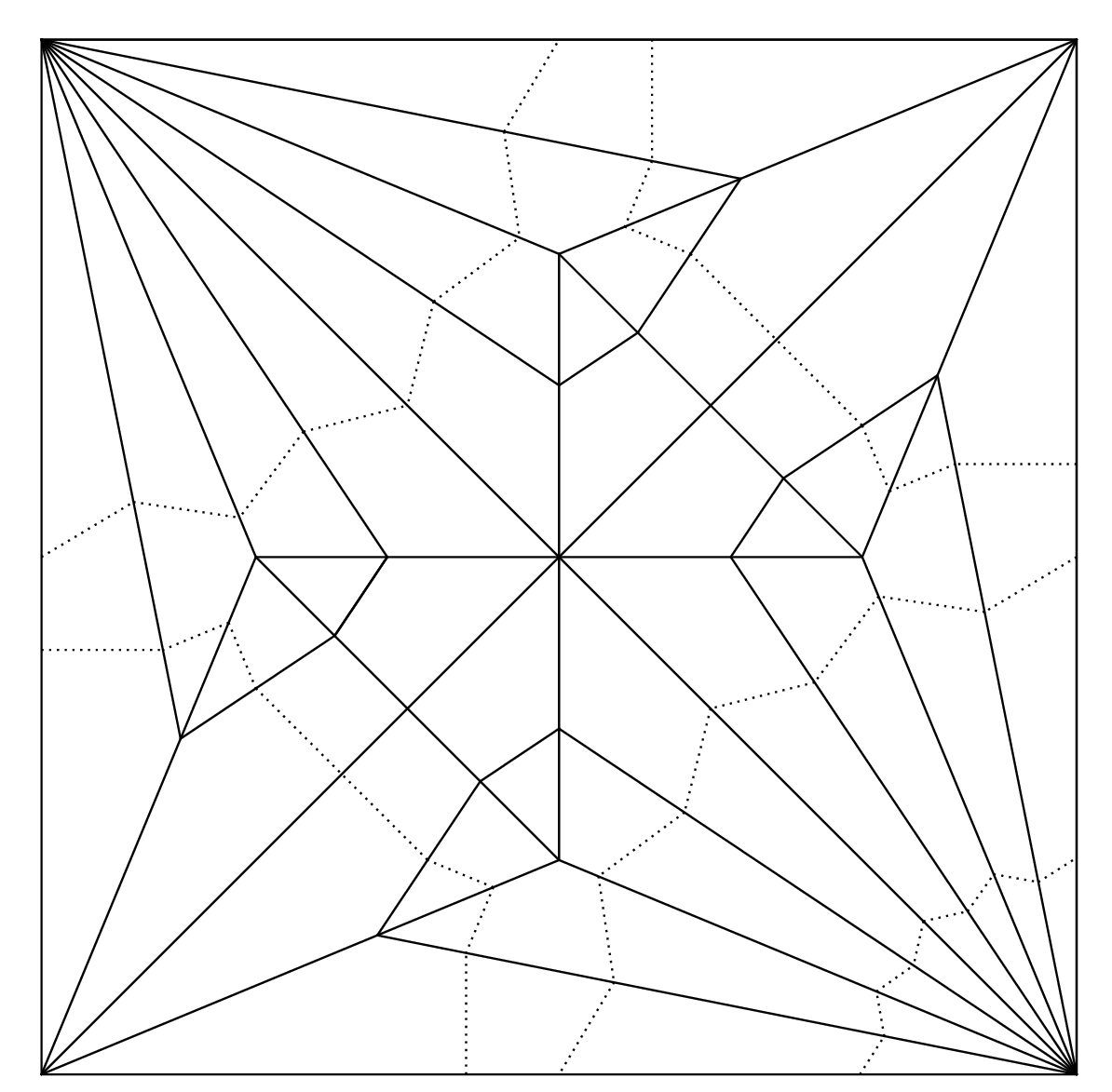}
\ \ \ \ \ \ \ \ 
\end{figure}

The base for the crane has four long ends and one short end.
Two ends are used for wings and the other two have to be thinned,
as one is used for the head and the other for the tail.
Thinning each of the long ends twice makes it possible to produce a base which can then be folded into the plane to obtain a figure with a larger perimeter.
\begin{center}
\begin{lpic}[t(0mm),b(0mm),r(0mm),l(0mm)]{pics/zagotovka(0.24)}
{\large
\lbl{230,170,-10;$\longrightarrow$}
\lbl{295,105,-20;$\longrightarrow$}
}
{\tiny
\lbl{119,119;$1$}
\lbl{134,110;$2$}
\lbl{147,110;$3$}
\lbl{160,110;$4$}
\lbl{170,110;$5$}
\lbl{180,110;$6$}
\lbl{190,110;$7$}
\lbl{200,110;$8$}

\lbl{117,100;$16$}
\lbl{134,100;$15$}
\lbl{147,100;$14$}
\lbl{160,100;$13$}
\lbl{170,100;$12$}
\lbl{180,100;$11$}
\lbl{190,100;$10$}
\lbl{200,100;$9$}

%\lbl{80;$1$}
\lbl{110,134;$2$}
\lbl{110,147;$3$}
\lbl{110,160;$4$}
\lbl{110,171;$5$}
\lbl{110,181;$6$}
\lbl{110,191;$7$}
\lbl{110,201;$8$}

\lbl{98,120;$16$}
\lbl{98,134;$15$}
\lbl{98,147;$14$}
\lbl{98,160;$13$}
\lbl{98,171;$12$}
\lbl{98,181;$11$}
\lbl{98,191;$10$}
\lbl{98,201;$9$}

\lbl{90,110;$17$}
\lbl{74,110;$18$}
\lbl{60,110;$19$}
\lbl{48,110;$20$}
\lbl{37,110;$21$}
\lbl{26.5,110;$22$}
\lbl{17,110;$23$}
\lbl{7,110;$24$}

\lbl{74,100;$31$}
\lbl{60,100;$30$}
\lbl{48,100;$29$}
\lbl{37,100;$28$}
\lbl{27,100;$27$}
\lbl{17,100;$26$}
\lbl{7,100;$25$}

\lbl{110,88;$17$}
\lbl{110,74;$18$}
\lbl{110,59;$19$}
\lbl{110,48;$20$}
\lbl{110,38;$21$}
\lbl{110,28;$22$}
\lbl{110,18;$23$}
\lbl{110,8;$24$}

\lbl{88,88;$32$}

\lbl{98,74;$31$}
\lbl{98,59;$30$}
\lbl{98,48;$29$}
\lbl{98,38;$28$}
\lbl{98,28;$27$}
\lbl{98,18;$26$}
\lbl{98,8;$25$}
}
\end{lpic}
\end{center}
In the picture above, you can see the net of folds, the base, and the base with opened out ends.
On the net of folds, each region is labeled with the number of the corresponding layer in the base.
The dashed lines are the folds which appear when the ends are opened out.
The perimeter increases by about $0.5\%$, and there are 80 layers in the end.
We do not know of a way to increase the perimeter with a smaller number of layers.

If $a$ is the side length of the original square, then it takes a bit less than $a$ to go around each of the four needles and it takes about $(\sqrt{2}-1)\cdot a$ to go around the short end, resulting in a longer perimeter.
Thinning the ends many times makes it possible to increase the perimeter by a value arbitrarily close to $(\sqrt{2}-1)\cdot a$.

\medskip

The following picture describes another way to increase the perimeter, based on an idea of Yashenko \cite{yashenko}.
It can be obtained by recursive application of one simple move.
Repeating this move sufficiently many times produces a figure with a longer perimeter.
Since each iteration adds two layers near the concave corner,
the total number of layers in this model is much larger than in the crane base.

\begin{center}
\begin{lpic}[t(0mm),b(0mm),r(0mm),l(0mm)]{pics/yaschenko(0.54)}
\end{lpic}
\end{center}

\parbf{The sea urchin and the comb.}
It turns out that the perimeter can be made arbitrarily large.
This can be seen in the origami model for a sea urchin constructed by Robert~Lang in 1987 \cite{lang}.
In 2004, a complete solution was discovered independently by Alexei Tarasov \cite{tarasov}.
Tarasov constructs a folding of a ``comb'', which is shown in the picture on the right.
More importantly, he proves that the comb can be folded in a true way,
in particular without starching and crooking the paper as is often done in origami.

\begin{center}
\ \ \ \ 
\includegraphics[scale=0.28]{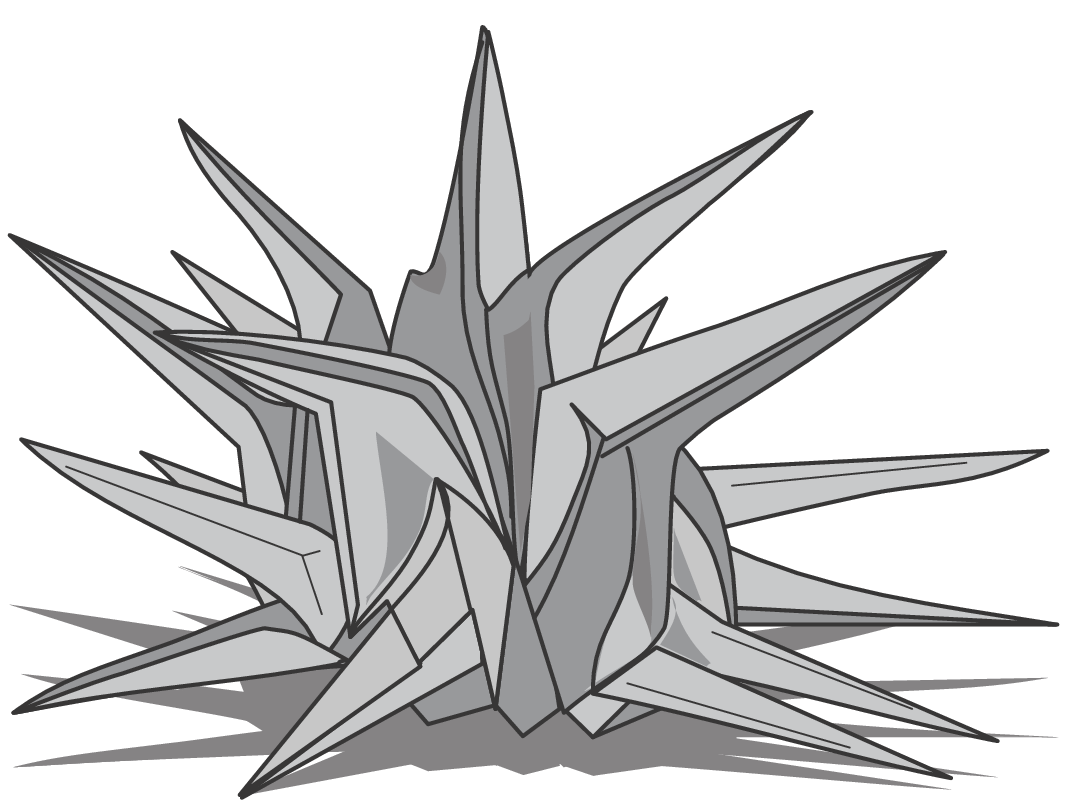}
\hfill
\includegraphics[scale=0.55]{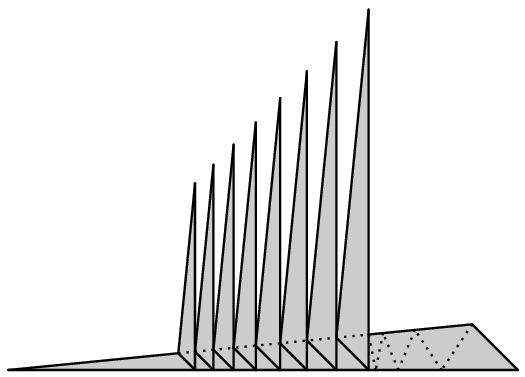}
\ \ \ \ 
\end{center}

If you read the original pdf-file on a computer, 
you can extract the following three movies which describe Tarasov's solution, 
\begin{itemize}
\item \textattachfile[color=0 0 1,author=Alexey Tarasov,mimetype=video/avi]{pics/rouble.avi}{The first movie} shows the folding of the whole comb.
\item \textattachfile[color=0 0 1,author=Alexey Tarasov,mimetype=video/avi]{pics/layer-up.avi}{The second movie} shows the folding of one of the needles separately.
\item \textattachfile[color=0 0 1,author=Alexey Tarasov,mimetype=video/avi]{pics/2layers.avi}{The third movie} shows the folding of two needles simultaneously.
\end{itemize}

\parbf{Foldings in 4-dimensional space.} 
One can define a ``folding'' as a piecewise distance-preserving map from the square to the plane.
These foldings are yet more general than those which appear above.
The following exercise 
%shows that it is not always possible to realize  such a map as a folding.
shows that it is not always possible to realize  such a map by folding a paper model.

\begin{wrapfigure}{r}{40mm}
\begin{lpic}[t(-7mm),b(-0mm),r(0mm),l(0mm)]{pics/3464(0.4)}
\end{lpic}
\end{wrapfigure}

\begin{thm}{Exercise}\label{pr:6-4-3-4}
Consider the portion of a regular tessellation of the square $\square$ as in the picture.

Show that there is a map $f\:\square\to\RR^2$ which is distance-preserving on each polygon in the tessellation
and such that it
only reverses the orientation of the gray polygons.%
\footnote{Less formally, you need to ``fold'' along each segment in this tessellation.}

Show that it is not possible to make a paper folding model for $f$.%
\footnote{More formally, we need to think of the plane as lying in $\RR^3$,
and we need to show that the map $f$ cannot be approximated by injective continuous maps $\square\to \RR^3$.}
\end{thm}

The obstructions described in the above exercise disappear in $\RR^4$;
i.e., one can regard piecewise distance-preserving maps as paper foldings in 4-dimensional space.
Moreover, one can actually fold the square in $\RR^4$, as prescribed by a given piecewise distance-preserving map.
By this we mean that one can construct a continuous one-parameter family of piecewise distance-preserving maps $f_t\:\square\to\RR^4$, $t\in[0,1]$, with a fixed triangulation, say $\mathcal{T}$,
such that
\begin{itemize}
\item $f_0$ is a distance-preserving map from $\square$ to the coordinate plane $\RR^2\times \{0\}$ in $\RR^4=\RR^2\times\RR^2$,
\item the map $f_1$ is our given piecewise distance-preserving map to the same coordinate plane,
\item the map $f_t$ is injective for any $t\ne1$.
\end{itemize}

The proof of the last statement is based on Exercise~\ref{pr:alexander}.
Let $a_1,\dots,a_k$ 
be the vertices of $\mathcal{T}$
and $b_1,\dots,b_k$ be the corresponding images under the piecewise distance-preserving map.
Set $f_t(a_i)\in\RR^2\times\RR^2$ to be 
$$f_t(a_i)= \left(\frac{a_i + b_i}{2} + 
\cos(\pi\cdot t)\cdot \frac{a_i - b_i}2,\  
\sin(\pi\cdot t)\cdot \frac{a_i - b_i}2\right);
$$
so $f_0(a_i)=(a_i,0)$ and $f_1(a_i)=(b_i,0)$ for any $i$.
We can extend $f_t$ linearly to each triangle of $\mathcal{T}$.
Direct calculations show that $\ell_{i,j}(t)=|f_t(a_i)-f_t(a_j)|$ is  monotonic in $t$;
in particular, if $|a_i-a_j|=|b_i-b_j|$,
then $\ell_{i,j}(t)$ is constant.
This proves that $f_t$ is piecewise distance-preserving.

Further, direct calculations 
show that for any $x,y\in\square$,
$$|f_t(x)-f_t(y)|^2=p-q\cdot\cos(\pi\cdot t)$$
for some constants $p$ and $q$.
Therefore, if $x\ne y$, then $|f_t(x)-f_t(y)|>0$ for any $t\ne1$.
In other words, $f_t$ is injective for any $t\ne1$.

\medskip
 
It follows that for paper folding in $\RR^4$,
the existence of perimeter-increasing folds follows from Brehm's theorem.
It is sufficient to construct a distance non-expanding map $f$ from the square to the plane
so that the perimeter of its image is sufficiently large.
Then, applying Brehm's theorem to a sufficiently dense finite set of points in the square,
we get a piecewise distance-preserving map $h$ which is arbitrarily close to $f$.
In particular, we can arrange it so that the perimeter of the image $h(\square)$
is still sufficiently large.

The needed map can be constructed as follows:
Fix a large $n$ and divide the square $\square$ into $n^2$ squares with side length $\tfrac{a}{n}$.
Let $d(x)$ denote the distance from a point $x\in\square$ to the boundary of the small square which contains $x$.
The function $d\:\square\to\RR$ takes values in $[0,\tfrac{a}{2{\cdot}n}]$.
Further, let us enumerate the squares by integers from $1$ to $n^2$.
Given $x\in \square$,
denote by $i(x)$ the smallest number assigned to a small square which contains $x$.

\begin{wrapfigure}{r}{40mm}
\begin{lpic}[t(-5mm),b(-0mm),r(0mm),l(0mm)]{pics/square(1)}
\lbl[b]{27,11.6;$\xto{\ f\ }$}
\end{lpic}
\end{wrapfigure}

Now for each $i\in\{1,\dots,n^2\}$ choose a unit vector $u_i\in\RR^2$ so that $u_i\not=u_j$ if $i\not=j$.
Consider the map $f\:\square\to\RR^2$ defined by 
$$f(x)=d(x)\cdot u_{i(x)}.$$
It is straightforward to check that the resulting map is distance non-expanding.
The image $f(\square)$ consists of $n^2$ segments of length $\tfrac{a}{2{\cdot}n}$ which start at the origin.
So the perimeter of $f(\square)$ is equal to 
$2{\cdot}n^2{\cdot} \tfrac{a}{2{\cdot} n}=a{\cdot} n$.
So by taking $n$ large, one can make the perimeter of the image $f(\square)$ arbitrarily large.

(The picture shows the case $n=4$.
In this case the perimeter of $f(\square)$ is $4\cdot a$,
which is the same as the perimeter of the original square.
However, for $n>4$, it gets larger.
When you calculate the perimeter of the degenerate figure, imagine going up and down each segment as you ``traverse the boundary'' and counting the length of each segment twice.)

\addtocontents{toc}{\protect\end{quote}}
\chapter*{Final remarks}\label{sec:pdp-comments}
\addcontentsline{toc}{chapter}{Final remarks}
\markboth{\MakeUppercase{Final remarks}}{}

\parbf{Zalgaller's folding theorem.}
Zalgaller's theorem holds in all dimensions: any $m$-dimensional polyhedral space $P$ admits a piecewise distance-preserving map to $\RR^m$.

In \cite{zalgaller-PL}, Zalgaller proved this statement for $m\le 4$.
The trick described in the ``proof with no cheating'' makes the proof work in all dimensions.
This trick first appeared in Krat's thesis \cite{krat}.

\parbf{Brehm's extension theorem.}
This was proved  by Brehm \cite{brehm}
and rediscovered independently many years later by Akopyan and Tarasov \cite{akopyan-tarasov}.
The proofs are based on the same idea.

Brehm's extension theorem holds in all dimensions
and it can be proved along the same lines.

\parbf{Kirszbraun theorem.}
This remarkable theorem was proved by Kirszbraun in his thesis, defended in 1930.
A few years later he published the result in \cite{kirszbraun}.
Independently the same result was reproved later by  Valentine \cite{valentine}.

The paper of Danzer, Gr\"unbaum and Klee \cite{DGK}, which is delightful to read, gives a proof of this theorem  based on Helly's theorem on the intersection of convex sets.

\parbf{Akopyan's approximation theorem.}
This theorem 
also admits direct generalizations to higher dimensions.
Moreover the condition that $f$ is piecewise linear is redundant.
This is because
any distance non-expanding map from a polyhedral space to $\RR^m$ can be approximated by piecewise linear distance non-expanding maps.

The 2-dimensional case was proved by Krat 
in her thesis \cite{krat}.
In~\cite{akopyan},
Akopyan noticed that Brehm's extension theorem simplifies the proof and also makes it possible 
to prove the higher dimensional case.

Much earlier, an analogous question was considered by Burago and Zalgaller.
They proved that any piecewise linear embedding 
of a $2$-dimensional polyhedral surface in $\RR^3$ 
can be approximated by a piecewise distance-preserving embedding; see \cite{burago-zalgaller-0,burago-zalgaller}.

\parbf{Rumpling the sphere.}
Theorem \ref{thm:S2->R2} 
admits the following generalization,
which can be proved along the same lines.

{\sloppy

\begin{thm}{Gromov's rumpling theorem}
Let $M$ be an $m$-dimensional Riemannian manifold.
Then any distance non-expanding map $f\:M\z\to \RR^m$
can be approximated by length-preserving maps.
More precisely, given $\eps>0$ there is a length-preserving map $f_\eps\:M\to\RR^m$
such that 
$$|f_\eps(x)-f(x)|<\eps$$
for any $x\in M$.
\end{thm}

}

This result is a partial case of Gromov's theorem in \cite[Section~2.4.11]{gromov}.
The proof presented here is based on the construction in \cite{petrunin-inverse};
Gromov's original proof is different.
The proof in \cite{petrunin-inverse}
also makes it possible to construct surprising examples of spaces which admit length-preserving maps to $\RR^m$, such as sub-Riemannian manifolds.

Gromov's theorem states that length-preserving maps
have no non-trivial global properties.
This is a typical ``local to global'' problem.
Here
the length-preserving property is ``local''
and the only ``global'' consequence is trivial:
it is the distance non-expanding property.

For such ``local to global'' problems the answer
``no non-trivial global properties''
 is the most common, 
but
it does not mean that it is easy to prove.
The so-called ``\textit{h}-principle'' (homotopy principle) provides machinery for proving such statements.
The \textit{h}-principle is not a theorem; it is a property which often holds for different geometric structures.
There are a few methods to prove the
\textit{h}-principle, including the one which is described in the proof of Theorem \ref{thm:S2->R2}%
\footnote{Usually, the \textit{h}-principle is formulated in terms of partial differential equations, 
but one may think of ``length-preserving maps'' as weak solutions of a particular partial differential equation.}.
Gromov's rumpling theorem is one of the simplest examples.
Other examples include 
\begin{itemize}
\item The cone eversion theorem, which states that there is a continuously varying one-parameter family of smooth functions $f_t(x,y)$, $t\in[0,1]$, without critical points in
the punctured plane $\RR^2\backslash\{0\}$, 
such that
$f_0(x,y)=-\sqrt{x^2+y^2}$ and $f_{1}(x,y)=\sqrt{x^2+y^2}$.
See \cite[Lecture 27]{TF} and read the whole book; it is nice.
\item The Nash--Kuiper theorem, which in particular implies the existence of $C^1$-smooth length-preserving maps $\SS^2\to\RR^3$ whose image has arbitrarily small diameter.
\item Smale's sphere eversion paradox, which states that there is a continuous one-parameter family of smooth immersions $f_t\:\SS^2\z\to \RR^3$, $t\in[0,1]$ such that $f_0\:\SS^2\z\to \RR^3$ is the standard inclusion and $f_1(x)=-f_0(x)$ for all $x\in \SS^2$.
\item Combining the techniques of Smale and Nash--Kuiper, 
one can make sphere eversions $f_t$ in the class of length-preserving $C^1$-smooth maps.
\end{itemize}
 
For further reading we suggest the comprehensive introduction to the \textit{h}-principle 
by Eliashberg and Mishachev \cite{eliashberg-mishachev}.

\parbf{Paper folding.} The aspects of paper folding related to geometric constructions are discussed in \cite{hull};
this paper is very entertaining.
Interesting aspects of paper folding in the 3-dimensional space are covered in \cite[Lecture 15]{TF}.

%\parbf{Acknowledgments.} 
%We would like to thank
%Arseniy Akopyan, 
%Robert Lang, 
%Alexei Tarasov
%for their help.
%Also we would like to thank all the students in our class
%for their participation and true interest.

\chapter*{Hints and solutions}
\addcontentsline{toc}{chapter}{Hints and solutions}
\refstepcounter{chapter}
\markboth{\MakeUppercase{Hints and solutions}}{}

\parbf{Exercise \ref{ex:IsometriesOfR2}; (a).}
Given an isometry $\iota: \RR^2 \to \RR^2$, let
$$F_{\iota} = \set{x \in \RR^2}{\iota(x) = x}$$ be the set of fixed points of $\iota$.
Show that if $x, y \in F_{\iota}$, then the line thru $x$ and $y$ is contained in $F_{\iota}$.
Conclude that for any isometry, $F_{\iota}$ is either empty, a point, a line, or all of $\RR^2$.

Given two isometries $\iota_1$ and $\iota_2$ that agree on three non-collinear points, use the above to argue that
$$\iota_1 \circ (\iota_2)^{-1} = \id_{\RR^2},$$ the identity map on $\RR^2$.

\parbf{(b).}
Notice that translations by a fixed vector, rotations about a point, and reflections across a line are all examples of isometries of $\RR^2$.
The required isometry can be constructed as the composition of a translation, followed by a rotation, and then (possibly) a reflection.

\parbf{Exercise \ref{LP=>short}. (\ref{LP=>short:a})} follows from the definition of length space \ref{def:length-space}.

\parbf{(\ref{LP=>short:b}).}
Suppose $f$ is distance non-expanding.
It follows from the definition of length that $\length (f\circ \alpha) \le \length \alpha$.

For the reverse inequality, let $\eps > 0$ be arbitrary and choose a partition $t_0 < t_1 < \dots < t_n$ such that
$$\length \alpha - \eps < \sum_{i=1}^n |\alpha(t_i) - \alpha(t_{i-1})|_X.$$
Let $\alpha_i = \alpha|_{[t_{i-1}, t_i]}$, which is the arc of $\alpha$ from $\alpha(t_{i-1})$ to $\alpha(t_i)$.
Then using the assumption on $f$, we have
\begin{align*}
\length \alpha - \eps &< \sum_{i=1}^n |\alpha(t_i) - \alpha(t_{i-1})|_X \le 
\\
&\le \sum_{i=1}^n \length(f \circ \alpha_i) =
\\
&=\length(f\circ \alpha).
\end{align*}
Since $\eps >0$ was arbitrary, this shows $\length \alpha \le \length (f \circ \alpha)$.

\parbf{Exercise \ref{ex:n=<m}.}
Fix an $m$-dimensional simplex $\Delta$ and a distance-preserving map $f: \Delta \to \RR^n$.

Note that line segments in $\Delta$ are mapped to line segments in $\RR^n$.
Moreover perpendicular line segments are mapped to perpendicular line segments.
Both statements follow since they can be formulated entirely in terms of the metrics on the spaces.

Notice that $m$ is the maximal number of mutually perpendicular line segments that can pass thru a point in $\Delta$, and $n$ is the maximal number of mutually perpendicular line segments that can pass through a point in $\RR^n$.
Hence the statement follows.

\parbf{Exercise \ref{pdp-for-tetrahedron}.}
It would be helpful to glue a paper model of $\partial\Delta$ and try to fold it onto the plane.

\begin{wrapfigure}{r}{35mm}
\begin{lpic}[t(-6mm),b(-2mm),r(0mm),l(0mm)]{pics/surface-of-tetrahedron(1)}
\lbl[br]{2,11;$a$}
\lbl[lb]{33,11;$b$}
\lbl[b]{25.5,27;$c$}
\lbl[t]{8,2.5;$d$}
\lbl[lb]{21.5,13;$a'$}
\lbl[brw]{12,13.5;$b'$}
\lbl[lw]{19.5,17;$x$}
\end{lpic}
\end{wrapfigure}

Let $a$, $b$, $c$ and $d$ be the vertices of~$\Delta$.
If $f\: \partial \Delta \to \RR^2$
is distance-preserving on each face, then it is also distance-preserving on the subset $\{a, b, c, d\}$.
Thus, $\Conv\{f(a), f(b), f(c), f(d)\}$ is an isometric copy of $\Delta$ in $\RR^2$, which is impossible by Exercise \ref{ex:n=<m}.

{\sloppy

The picture shows a triangulation for a piecewise distance-preserving map $f\: \partial \Delta \to \RR^2$.
The boundary $\partial \Delta$ is subdivided into 10 triangles.
The points $a'$ and $b'$ are tangent points on the faces $\triangle bcd$ and $\triangle acd$
to the sphere inscribed in $\Delta$.
In particular $\triangle a'cd \cong \triangle b'cd$.

}

If $c'$ and $d'$ denote the tangent points on the faces $\triangle abd$ and $\triangle abc$ to the inscribed sphere, then we also have 
$\triangle a'bc \cong \triangle d'bc$, 
$\triangle a'bd \cong \triangle c'bd$
and so on, with 6 pairs of congruent triangles altogether.
This makes it possible to fold it onto the plane so that corresponding points in these pairs will coincide.
(In particular, all the points $a'$, $b'$, $c'$ and $d'$ will be mapped to one point.)

\parbf{Exercise \ref{ex:acute-triangulation}.}
The solution should be clear from the picture.

\parit{Comment.}
As shown in \cite{saraf}, any 2-dimensional polyhedral space admits a triangulation in which all of the triangles are acute.

\begin{wrapfigure}{r}{63mm}
\begin{lpic}[t(-0mm),b(-2mm),r(-0mm),l(-2mm)]{pics/acute-triangulation(1)}
\end{lpic}
\end{wrapfigure}

One may call a simplex \emph{acute}\index{acute simplex}
if it contains its own circumcenter
--- 
this provides a natural generalization of an acute triangle to higher dimensions.
The existence of acute triangulations in higher dimensions seems to be unlikely, but as far as we can see, nothing is known about these triangulations.

\parbf{Exercise \ref{ex:voronoi-in-star}.}
First note that if $x\in V_i$, then any geodesic $[w_i,x]$ lies in $V_i$.
Indeed if $y\in [w_i,x]$, then for any $j$ we have
\begin{align*}
|w_i-y|
&=|w_i-x|-|y-x|\le
\\
&\le |w_j-x|-|y-x|\le
\\
&\le 
|w_j-y|.
\end{align*}
This shows $y\in V_i$.

Therefore if $V_i\not\subset S(w_i)$
then there is a triangle $\Delta$ of the triangulation such that $w_i\in \Delta$
and $V_i$ contains a point $x$ on a side $E$ of $\Delta$
which does not contain $w_i$.

Choose $j$ so that $w_j\in E$ and $|w_j-x|$ is minimal.
Note that $|w_i-x|\le |w_j-x|<\eps$.
Since $\eps< \tfrac{\ell\cdot \alpha}{100}$,
there is a vertex, say $v$, of $\Delta$ such that
$|w_i-v|<\tfrac\ell2$.
Therefore there is $w_n\in E$ such that $|v-w_i|=|v-w_n|$.
Finally note that 
$$|w_i-x|>\big| |v-x|-|v-w_i| \big|=|w_n-x|;$$
i.e. $x\notin V_i$, a contradiction.

\parbf{Exercise \ref{ex:zalgalle+embedding}.}
In other words, we need to show that there is an injective piecewise distance-preserving map from the $2$-dimensional polyhedral space $P$
into a Euclidean space of sufficiently large dimension.

Fix a triangulation $f\:|K|\to P$ of $P$, where $K$ is a simplicial complex in some $\RR^n$.
We can assume that the triangulation is linear,
that is, $f$ sends a point in $K$ to the corresponding point in $P$ with the same barycentric coordinates.

{\sloppy

Equip the underlying set $|K|$ with the induced length metric from~$\RR^n$.
The complex $K$ can be rescaled to ensure that the map $f\:|K|\z\to P$ is distance-expanding,
meaning there is $\lambda < 1$ such that if $x, y\in |K|$, then
\[|x-y|_{|K|}\le \lambda\cdot|x'-y'|_{P},\]
where $x' = f(x)$ and $y' = f(y)$.

}

Show that there is a unique length metric $\rho$ on $P$
such that for any two points $x,y$ in one simplex $\triangle$ of $K$
we have \[\rho(x',y')=\sqrt{|x'-y'|_P^2-|x-y|_{|K|}^2},\]
where again $x'=f(x)$ and $y'=f(y)$.

Note that $(P,\rho)$ is a polyhedral space.
Applying Zalgaller's theorem, we get a piecewise distance-preserving map $h\:(P,\rho)\to \RR^2$.
Note that the map $P\to \RR^2\times\RR^n=\RR^{n+2}$ defined by 
\[x\mapsto (h(x),f^{-1}(x))\] is injective and piecewise distance-preserving.
Hence the statement follows.

\parbf{Exercise \ref{ex:black-and-white}.}
First note that the statement is trivial
if the triangulation has only one interior vertex.

Order the triangles of the triangulation 
in such a way that each triangle intersects 
the union of the previous triangles on one or two sides.
(This might look obviously possible, 
but try to give a proof.
In more sophisticated language,
our triangulation is \emph{shellable}; you can search for this term on the web.)

Let us fix the map on the first triangle
and extend it to the subsequent triangles in order by ``folding'' (i.e., reflecting) along the sides where the color changes.
Note that on each step, if the map exists, then it has to be unique.

The existence might only fail if the new triangle 
has two common sides with the old triangles.
In this case one only has to check a neighborhood of the common vertex of these two sides.
In this way we reduce to the case of the triangulation with only one interior vertex.

\begin{wrapfigure}{r}{48mm}
\begin{lpic}[t(-0mm),b(-5mm),r(0mm),l(0mm)]{pics/Q12(1)}
\lbl[tr]{2,1;$a_1$}
\lbl[tl]{43.5,1;$x_1$}
\lbl[bl]{44,22;$x_2$}
\lbl[b]{28,37,-10;$\dots$}
\lbl[br]{8,33;$x_k$}
\lbl[tl]{27,21;$y_k$}
\lbl[br]{22,38;$\ell$}
\lbl{10,25;$Q_2$}
\lbl{30,10;$Q_1$}
\end{lpic}
\end{wrapfigure}

\parbf{Exercise \ref{ex:triangle-reflect}.} 
Without loss of generality, we may assume that
$a_1=b_1$, 
$y_1=x_1$, 
and both points $x_k$ and $y_k$ lie on the same side of the line $a_1x_1$.
These can be arranged by applying a translation, followed by a rotation, and then possibly a reflection.
The case where $|y_1-y_k| \z= |x_1-x_k|$ is handled by Exercise \ref{ex:IsometriesOfR2Existence}.
So assume that
$|y_1-y_k|\z<|x_1-x_k|$.

Let $\ell$ be the bisector of the angle $\angle x_k a_1 y_k$.
Note that the reflection of $x_k$ in $\ell$ is $y_k$.
Since 
\[|y_1-y_k|<|x_1-x_k|,
\]
 $\ell$ cuts the polygon $Q$
into two parts, say $Q_1\ni x_1$ and $Q_2\ni x_k$.

Let $f(x)=x$ if $x\in Q_1$, and let $f(x)$ be the reflection of $x$ in $\ell$ if $x\in Q_2$.

\begin{wrapfigure}{r}{65mm}
\begin{lpic}[t(-1mm),b(-0mm),r(0mm),l(0mm)]{pics/bigger-area(1)}
\lbl[b]{13,-1;$a_1$}
\lbl[b]{3,23;$a_2$}
\lbl[b]{33,-1;$a_3$}
\lbl[b]{41.5,-1;$b_1$}
\lbl[b]{44,24;$b_2$}
\lbl[b]{61.5,-1;$b_3$}
\lbl[]{16,7;$A$}
\lbl[]{48,8;$B$}
\end{lpic}
\end{wrapfigure}

\parbf{Exercise \ref{pr:perimeter}.}%
\footnote{The first inequality follows easily from two famous theorems in discrete geometry:
one is Alexander's theorem \cite{alexander}
and the other is the 
Kneser--Poulsen conjecture, which was solved in the 2-dimensional case by Bezdek and Connelly in \cite{bezdek-connelly}.
(The reduction to each of these theorems takes one line, 
but it might not be completely evident.)
We encourage you to read these papers; they are totally beautiful.
You will be surprised to learn that to solve this 2-dimensional problem, it is convenient to work in 4-dimensional space.
The reason is explained in Exercise~\ref{pr:alexander}.}
The second inequality does not hold in general.
This can be seen in the picture.
Below we give a proof of the first inequality from \cite{petrunin-ruble}.

Given a finite collection of points $a_1,\dots,a_n$,
set 
$$\ell(a_1,\dots, a_n)=\per[\Conv\{a_1,\dots, a_n\}].$$
Then $\per A\ge\per B$ can be written as
$$\ell(a_1,\dots, a_n)\ge \ell(b_1,\dots, b_n).$$
Applying Brehm's extension theorem (\ref{thm:brehm}),
we get a piecewise distance-preserving map $f\:A\to\RR^2$
such that $f(a_i)=b_i$ for each $i$.

Assume to the contrary that
$$\ell(a_1,\dots, a_n)< \ell(b_1,\dots, b_n).\eqlbl{eq:conta-a-b}$$
We can assume that 
$\{a_1,\dots,a_n\}$ and $f$ are chosen in such a way that 
\begin{clm}{}
\label{min-n} The number $n$ is the minimal value for which \ref{eq:conta-a-b} can hold.
\end{clm}
and
\begin{clm}{}\label{max-l} If $x_1,\dots, x_n\in A$
and $y_i=f(x_i)$ then
$$\ell(y_1,\dots,y_n)-\ell(x_1,\dots,x_n)\le\ell(b_1,\dots, b_n)- \ell(a_1,\dots, a_n).$$

\end{clm}
\noi
To meet Condition~\ref{max-l}, 
one has to replace $\{a_1,\dots,a_n\}$ with
an $n$-point subset $\{a_1',\dots, a_n'\}$ of $A$ 
that maximizes
$$\ell(f(a_1'),\dots,f(a_n')) - \ell(a_1',\dots, a_n').$$
This is possible since this difference is continuous and $A$ is closed and bounded.

Note that all $b_i$ are distinct vertices of $B$.
Indeed, if $b_n$ lies inside or on a side of $B$, then% or $b_n=b_i$ for some $i\not=n$.
\begin{align*}
\ell(b_1,\dots,b_{n-1})&=\ell(b_1,\dots,b_{n-1},b_{n}),
\\
\ell(a_1,\dots,a_{n-1})&\le \ell(a_1,\dots,a_{n-1},a_{n}).
\end{align*}
This contradicts Condition~\ref{min-n}.

\begin{wrapfigure}{r}{41mm}
\begin{lpic}[t(-0mm),b(-0mm),r(0mm),l(0mm)]{pics/obkhvat-shorter(1)}
\lbl[tr]{2,15;$b_1$}
\lbl[b]{21,39;$b_2$}
\lbl[lt]{33,2;$b_3$}
\lbl[r]{2,32;$b_4$}
\lbl[t]{15,2;$b_5$}
\lbl[bl]{31,36;$b_6$}
\lbl[l]{38,20;$b_7$}
\end{lpic}
\end{wrapfigure}

By $\mangle a_i$ and $\mangle b_i$, 
we will denote the angles 
of $A$ and $B$ at $a_i$ and $b_i$ respectively.
If $a_i$ lies inside or on a side of $A$ we set $\mangle a_i=\pi$.
Let us show that
\[\mangle b_i\le \mangle a_i.
\eqlbl{eq:a<b}\] 
If we move $a_i$ with unit speed inside $A$ along the angle bisector at $a_i$,
then the value $\ell(a_1,\dots,a_{n})$ 
decreases at a rate of $2{\cdot}\cos\tfrac{\mangle a_i}2$.
The point $b_i=f(a_i)$ will also move with unit speed.
One can show that the value $\ell(b_1,\dots,b_{n})$
cannot decrease at a rate greater than $2{\cdot}\cos\tfrac{\mangle b_i}2$.
By Condition~\ref{max-l}, 
the difference 
$$\ell(b_1,\dots,b_{n})-\ell(a_1,\dots,a_{n})$$
cannot increase.
Therefore
$2{\cdot}\cos\tfrac{\mangle b_i}2
\ge
2{\cdot}\cos\tfrac{\mangle a_i}2$,
hence \ref{eq:a<b}.

Applying the theorem about the sum of the angles of an $n$-gon to \ref{eq:a<b},
we get that $A$ is an $n$-gon with vertices $\{a_1,\dots,a_n\}$.
(We also get $\mangle b_i=\mangle a_i$, 
but we will not need it.)

By relabeling if necessary, we can assume that the points $a_i$ are labeled in the cyclic order in which they appear on the boundary of $A$.
Note that the same may not be true for the $b_i$.
In this case
$$\ell(a_1,\dots,a_{n})=|a_1-a_2|+\dots+|a_{n-1}-a_n|+|a_n-a_1|.$$
Since $|b_i-b_j|\le |a_i-a_j|$ for all $i$ and $j$,
we get
$$|b_1-b_2|+\dots+|b_n-b_1|\le|a_1-a_2|+\dots+|a_n-a_1|.$$
Finally, note that 
$$\ell(b_1,\dots,b_{n})\le|b_1-b_2|+\cdots+|b_n-b_1|.$$
The idea of the proof should be evident from the picture.

Therefore
$$\ell(b_1,\dots,b_{n})\le\ell(a_1,\dots,a_{n}),$$ 
a contradiction.

\parbf{Exercise \ref{pr:alexander}.}
Here is an example of such curves:
$$
\alpha_i(t) = \left(\frac{a_i + b_i}{2} + 
\cos(\pi\cdot t)\cdot \frac{a_i - b_i}2,\
\sin(\pi\cdot t)\cdot \frac{a_i - b_i}2\right).
$$
It is straightforward to check that
$\ell_{i,j}$ are monotonic.

\parbf{Exercise \ref{pr:brehm}.}
Let $A = \Conv Q$ and let $f: A \to \RR^2$ be the map produced by Brehm's theorem.
According to Exercise~\ref{problem2}, 
there is a distance non-expanding map $h\:\RR^2\to A$
such that $h(a)=a$ for any $a\in A$.
Taking the composition $F = f\circ h$, we get the needed distance non-expanding map $F: \RR^2\to \RR^2$.

\parbf{Exercise \ref{ex:PDPisPL}.}
First show that if $\Delta$ is a Euclidean $2$-simplex, then $f: \Delta \to \RR^2$ is linear if and only if the restriction of $f$ to any line segment in $\Delta$ is linear.
Use this to show that a distance-preserving map is linear.

\parbf{Exercise \ref{ex:akopyan-brehm}.}
Choose a sufficiently fine triangulation of $P$, 
say the diameter of each triangle is less than $\eps$.
If $\{a_1,\dots,a_n\}$ is the set of vertices of this triangulation,
take $b_i=h(a_i)$ and apply Brehm's extension theorem.
We obtain a map $f\:P\to\RR^2$ which coincides with $h$ on the set $\{a_1,\dots,a_n\}$.

Then the statement follows since the triangulation is fine 
and both $h$ and $f$ are distance non-expanding.

\parbf{Exercise \ref{ex:tripod}.}
Assume the contrary, and let $f$ be a piecewise distance-preserving map
from the tripod to the plane which fixes $a$, $b$ and $c$.
Since $f$ is distance non-expanding, we also get that $f(o)=o$.

\begin{wrapfigure}{r}{35mm}
\begin{lpic}[t(-2mm),b(-0mm),r(0mm),l(0mm)]{pics/abcox(1)}
\lbl[l]{34,22;$a$}
\lbl[br]{1,43;$b$}
\lbl[tr]{1,1;$c$}
\lbl[lb]{12,24;$o$}
\lbl[t]{15,7;$x'$}
\end{lpic}
\end{wrapfigure}

Take $x$ on the edge $[o,p]$ and set $x'\z=f(x)$.
Note that 
\begin{align*}
|x'-a|&\le |x-a|,\\
|x'-b|&\le |x-b|,\\
|x'-c|&\le |x-c|.
\end{align*}

Moreover, if we assume that $x$ is sufficiently close to $o$,
then 
\[|x'-o|=|x-o|.\]

It follows that
\begin{align*}
\measuredangle aox'&\le \measuredangle aox,
\\
\measuredangle box'&\le \measuredangle box,
\\
\measuredangle cox'&\le \measuredangle cox.
\end{align*}
However
\[\measuredangle aox=\measuredangle box=\measuredangle cox=\tfrac\pi2,\]
and therefore
\[\measuredangle aox',\measuredangle box',
\measuredangle cox'\le \tfrac\pi2.\]

On the other hand, it is clear that for any point $x'\ne o$
in the plane of $\triangle abc$,
at least one of the values $\measuredangle aox',\measuredangle box',
\measuredangle cox'$ exceeds $\tfrac\pi2$,
a contradiction.

\parbf{Exercise \ref{problem2}.}\label{sol-problem2}
The set $K$ is bounded and closed, so by the Extreme Value Theorem there is a point $\bar x\in K$ 
which minimizes the distance $|\bar x-x|$.

Assume there are two distinct points of minimal distance, say $\bar x$ and~$\bar x'$.
From convexity, their midpoint $z=\tfrac{\bar x+\bar x'}{2}$ lies in $K$.
Clearly 
$$|x-z|<|x-\bar x|=|x-\bar x'|,$$
a contradiction.

\begin{wrapfigure}{r}{45mm}
\begin{lpic}[t(-5mm),b(-3mm),r(0mm),l(0mm)]{pics/hilbert(1)}
\lbl[l]{43,18;$x$}
\lbl[br]{8,41;$y$}
\lbl[tl]{33,14;$\bar x$}
\lbl[br]{10,26;$\bar y$}
\lbl[br]{41,35;$\Pi_x$}
\lbl[tl]{19,42;$\Pi_y$}
\end{lpic}
\end{wrapfigure}

It remains to show that 
$$|\bar x-\bar y|\le|x-y|
\eqlbl{eq:|xy|=<|xy|}$$ 
for any $x,y\in\RR^3$.
We can assume that $\bar x\ne \bar y$;
otherwise, there is nothing to prove.

Consider the two planes $\Pi_x$ and $\Pi_y$ which pass thru $\bar x$ and $\bar y$ and are perpendicular to the line segment $[\bar x,\bar y]$.
Note that $x$ and $\bar y$ lie on opposite sides of $\Pi_x$;
otherwise, there would be a point on $[\bar x,\bar y]$ which is closer to $x$ than $\bar x$.
This is not possible since $[\bar x,\bar y]\subset K$.
In the same way we see that $y$ and $\bar x$ lie on opposite sides of $\Pi_y$.

Therefore the segment $[x,y]$ has to intersect both planes $\Pi_x$ and $\Pi_y$.
It remains to note that the distance from any point on $\Pi_x$ to any other point on $\Pi_y$ is at least $|\bar x-\bar y|$.
Hence \ref{eq:|xy|=<|xy|} follows.

\parbf{Exercise \ref{ex:limit-above}.}
Set $y_{n,k}=z_k$ for $n> k$.

Note that if $n>k$ then $z_k\in K_n$
and therefore $z_k=\phi_n(z_k)$.
I.e., the identity 
$$y_{n,k}=\phi_{n}(y_{n+1,k})$$ 
still holds for all $n$ and $k$.

Fix $k$ and $m$.
According to Exercise \ref{problem2},
the sequence 
$$\ell_n=|y_{n,k}-y_{n,m}|$$
is nondecreasing.
Since $\ell_n=|z_k-z_m|$ for $n>\max\{m,k\}$,
we get
$$|y_{n,k}-y_{n,m}|\le |z_k-z_m|$$
for all $n$.
It follows that for any fixed $n$,
the sequence
$\left(y_{n,m}\right)_{m=n}^\infty$ is a Cauchy sequence of points in $P_n$, 
and thus has a limit $x_n \in P_n$.
Since $\phi_n$ is continuous, we get 
$$\phi_n(x_{n+1}) \z= \phi_n\left(\lim_{m\to\infty} y_{n+1,m} \right) \z= \lim_{m\to\infty} y_{n,m} \z= x_n$$
for all $n$.

Let $\eps > 0$; choose $N$ sufficiently large so that $|z_m - z_n| < \tfrac{\eps}{3}$ for all $m,n \ge N$, and $|x - z_n| < \tfrac{\eps}{3}$ for all $n \ge N$.
Let $n \ge N$ and choose $m \ge N$ such that $|x_n - y_{n,m}| < \tfrac{\eps}{3}$.
From the above argument,
we have
\begin{align*}
|z_n - y_{n,m}| &= |y_{n,n} - y_{n,m}| \le
\\
&\le|z_n - z_m| <
\\
&<\tfrac{\eps}{3}.
\end{align*}
Therefore
$$
|x - x_n|\le |x - z_n| + |z_n - y_{n,m}| + |y_{n,m} - x_n|< \eps.
$$
It follows that $x_n \to x$ as $n\to\infty$.

\parbf{Exercise \ref{pr:6-4-3-4}.}
The existence of $f$ also follows from Exercise~\ref{ex:black-and-white}.

\begin{center}
%\begin{wrapfigure}{r}{62mm}
\begin{lpic}[t(-0mm),b(0mm),r(0mm),l(0mm)]{pics/triangle-3434(1)}
\lbl{35,10;$\longrightarrow$}
\end{lpic}
%\end{wrapfigure}
\end{center}

Note that the tessellation can be obtained by recursively reflecting the triangle in the picture across its sides and the sides of the resulting triangles.
It should be easy to construct the map on one triangle.
Note that this map sends each side of the triangle 
to a line which forms a side of the smaller equilateral triangle;
see the bold lines on the right side of the diagram.
Visualize this map and extend it to the whole plane by recursively applying reflections in the sides of these smaller triangles.

\begin{wrapfigure}{r}{32mm}
\begin{lpic}[t(0mm),b(0mm),r(0mm),l(0mm)]{pics/3-squares(1)}
\end{lpic}
\end{wrapfigure}

Notice that it is already impossible to fold the 9-gon in the picture along all three sides of the triangle.
Indeed, after folding, two of the squares have to lie on the same side
of the triangle.
Then the square which lies between the triangle and the other square
has to pass thru a side of the triangle.

\newpage
\input{pdp-minimum.ind}

\end{document}